%% file: gps-main.tex
\RequirePackage[l2tabu, orthodox]{nag} %

\documentclass[onefignum,onetabnum]{siamart190516}

\input{gps-shared}
\setlist[enumerate]{leftmargin=.7in}

\hypersetup{
  pdftitle={Gaussian Process Subspace Regression for Model Reduction},
  pdfauthor={Ruda Zhang, Simon Mak, and David Dunson},
  pdfprintscaling=None,
}

\begin{document}

\maketitle

\begin{abstract}
  Subspace-valued functions arise in a wide range of problems,
  including parametric reduced order modeling (PROM).
  In PROM, each parameter point can be associated with a subspace,
  which is used for Petrov-Galerkin projections of large system matrices.
  Previous efforts to approximate such functions use interpolations on manifolds,
  which can be inaccurate and slow.
  To tackle this, we propose a novel Bayesian nonparametric model for subspace prediction:
  the Gaussian Process Subspace regression (GPS) model.
  This method is extrinsic and intrinsic at the same time:
  with multivariate Gaussian distributions on the Euclidean space,
  it induces a joint probability model on the Grassmann manifold,
  the set of fixed-dimensional subspaces.
  The GPS adopts a simple yet general correlation structure,
  and a principled approach for model selection.
  Its predictive distribution admits an analytical form,
  which allows for efficient subspace prediction over the parameter space.
  For PROM, the GPS provides a probabilistic prediction at a new parameter point
  that retains the accuracy of local reduced models,
  at a computational complexity that does not depend on system dimension,
  and thus is suitable for online computation.
  We give four numerical examples to compare our method to subspace interpolation,
  as well as two methods that interpolate local reduced models.
  Overall, GPS is the most data efficient,
  more computationally efficient than subspace interpolation,
  and gives smooth predictions with uncertainty quantification.
\end{abstract}


\section{Introduction}

In this paper we propose a method to solve the following formal problem.
Consider a subspace-valued mapping $f: \Theta \mapsto G_{k,n}$
from a parameter space $\Theta \subset \mathbb{R}^d$ to the Grassmann manifold $G_{k,n}$,
which is the set of all $k$-dimensional subspaces of the Euclidean space $\mathbb{R}^n$.
Given function evaluations at $l$ points,
$(\boldsymbol{\theta}_i, \mathfrak{X}_i = f(\boldsymbol{\theta}_i ))_{i=1}^l$,
construct a probabilistic surrogate model $g$ such that $g(\boldsymbol{\theta}_*)$
is a probability distribution on $G_{k,n}$ concentrated near $f(\boldsymbol{\theta}_*)$
for any point $\boldsymbol{\theta}_* \in \Theta$.

\subsection{Motivation}

Many phenomena in science and engineering
can be described by systems of partial differential equations (PDEs).
For accurate analysis and prediction,
these mathematical models usually need to be discretized and simulated numerically.
This has lead to the development of computational science and engineering,
with wide-ranging applications such as aeroelastic systems \cite{Amsallem2008},
structural systems \cite{Panzer2010,Amsallem2011},
turbomachinery \cite{Mak2018}, ocean modeling \cite{ZhangRD2020EnvEcon},
and biomedicine \cite{chen2020function}.

Yet, high-fidelity models must resolve multiple physics,
multiple scales, complex geometry, and stochasticity.
This leads to large-scale dynamical systems that incur major computational costs,
especially when they need to be solved repeatedly.
Other applications require real-time or embedded computing based on limited computational resources.
In both cases, one needs to reduce the cost of solving large systems of differential equations.
Reduced order modeling (ROM) approximates the full model with a reduced order model,
which is a much smaller system of differential equations
that takes significantly less time and storage to simulate.
ROM often provides a speedup of several orders of magnitude,
and has been used in many types of problems in scientific computing \cite{Benner2015}.

Many ROM methods have been developed,
which can be roughly categorized into two types:
time domain methods such as proper orthogonal decomposition (POD) \cite{Lumley1967},
dynamic mode decomposition (DMD) \cite{Schmid2008}, and
discrete empirical interpolation method (DEIM) \cite{Chaturantabut2010};
and frequency domain methods such as balanced truncation \cite{Moore1981}
and rational interpolation \cite{Gugercin2008,Antoulas2020,Hokanson2020}.
Most of these ROM methods can be formulated as Petrov-Galerkin projections,
which projects the model state space onto a low-dimensional subspace.
Such a low-dimensional subspace is called a reduced subspace,
and a basis of the subspace is called a reduced basis.

In many use cases, the full model itself depends on some parameters,
to allow variations in material, geometry, loading, initial conditions, or boundary conditions.
However, the accuracy of reduced models often declines quickly as parameters change,
so we want to develop a reduced model that is also a function of the parameters.
This is called parametric reduced order modeling (PROM),
which is useful for design, control, optimization, uncertainty quantification, and inverse problems.
For a comprehensive review of projection-based PROM methods, see \cite{Benner2015}.

\subsection{Previous methods}
\label{sub:literature}

One approach is to consider a projection-based ROM method
as a mapping that associates each parameter point with a reduced subspace.
Given reduced subspaces at a sample of the parameter space,
one may approximate this subspace-valued mapping
and predict the reduced subspaces at other parameter points.
Compared with using one reduced subspace for the entire parameter space,
this keeps the reduced model small and often more reliable \cite{Amsallem2008}.

A natural idea to solve this problem
is to interpolate among local reduced subspaces as a deterministic function of the parameters,
using traditional interpolation methods.
However, this is infeasible since the Grassmann manifold is not a vector space
and linear combinations are undefined.
To circumvent this difficulty, \cite{Amsallem2008} proposed a method
that takes the interpolation to tangent spaces of the Grassmann manifold, which are vector spaces.
It comes in three steps.
Given a target parameter point, it chooses a few nearby parameter points
and maps the associated reduced subspaces to the tangent space of one of them
via the Riemannian logarithm.
Then the tangent vectors are interpolated as a function of the parameters,
using any traditional interpolation method.
Finally, the interpolated tangent vector is mapped back to the Grassmann manifold
via the Riemannian exponential, which gives the predicted subspace.
We will refer to this method as \textit{subspace interpolation} in this paper.
Some later developments include
accelerating computation for special types of interpolation schemes \cite{Son2013},
and adaptation for complex-valued data \cite{Zimmermann2014}.

Subspace interpolation has seen great success in PROM due to its accuracy and flexibility.
But its computational cost generally scales with the size of the full model,
which limits its use in large-scale problems.
To avoid this limitation,
\cite{Panzer2010} proposed a method that directly interpolates the reduced models:
it first applies a congruence transformation to the reduced models,
and then interpolates the model matrices element-wise.
We will refer to this approach as \textit{matrix interpolation}.
Influenced by this work, \cite{Amsallem2011} proposed a method that 
interpolates the transformed matrices on a relevant matrix manifold,
e.g. the general linear group, in a procedure analogous to subspace interpolation.
We will refer to this approach as \textit{manifold interpolation}.
Since their prediction costs do not scale with the full model,
these methods are considered as suitable for online computation.
However, all three interpolation methods lack a clear rule in selecting
the reference point, other interpolation points, and the interpolation scheme.
This often leads to model misspecification which undermines accuracy.
Moreover, subspace and manifold interpolation are extrinsic to the underlying manifolds,
where distortion is another source of error.

Another type of method uses the Riemannian center of mass of weighted data points.
The global or local Riemannian center of mass is the set of global or local minimizers
of the sum of weighted squared Riemannian distances \cite{Afsari2011}.
As before, the parameter-dependent weights can use any interpolation scheme
such as splines \cite{Grohs2013} or Lagrange polynomials \cite{Sander2016},
both of which were introduced in the context of geodesic finite elements.
Similarly, within the statistics community,
\cite{Petersen2019} proposed global and local regression models
where predictors are in a Euclidean space and random responses in a metric space.
These methods are intrinsic, i.e. involving operations entirely on the manifold,
so they avoid the limitations of mapping to a tangent space.
However, their computation requires iterative algorithms for Riemannian optimization,
and only local minimizers can be found.
So far their use has been limited to low-dimensional manifolds,
and we are not aware of their application in PROM.

Zimmermann \cite{Zimmermann2019} reviewed interpolation methods on the Grassmann manifold
and other matrix manifolds that arise in model reduction.
More recently, he introduced Hermite interpolation of
parameterized curves on Riemannian manifolds \cite{Zimmermann2020}.
All these methods are deterministic,
while probabilistic methods for subspace approximation have not been explored in the literature.

\subsection{Contribution}

We propose a new Gaussian process (GP) model for the approximation of subspace-valued functions,
which we call the \textit{GPS} model.
Instead of using differential geometric structures of the Grassmann manifold
as in \cite{Amsallem2008},
the GPS uses matrix-variate Gaussian distributions on the Euclidean space
to induce a probability model on the Grassmann manifold.
Our method therefore yields a probabilistic prediction of the subspace response,
with intrinsic characterization of its predictive mean and uncertainty.
Specifically, the mean prediction is a $k$-subspace of the span of the local reduced bases combined,
and the latter also covers most of the predictive uncertainty.
This GP model is flexible and yet well-guided:
it can be used with any correlation function on the parameter space,
and the function form and hyperparameters can be optimized via specific model selection criteria.

The main advantages of our method are summarized as follows.
(1) \textit{Data efficient}: accurate prediction requires only a small sample size $l$,
even when subspace dimension $k$ and parameter dimension $d$ are large.
(2) \textit{Computationally efficient}: its prediction cost does not depend on system dimension $n$,
and thus it is suitable for large-scale problems and online computation.
(3) \textit{Flexible}: It is a Bayesian nonparametric model
that is robust against model misspecification.
(4) It provides \textit{uncertainty quantification}, which gives confidence on a predicted subspace.

In our observation, GPS is much more accurate than subspace interpolation,
which is in turn much more accurate than manifold and matrix interpolation.
Such data efficiency can be attributed to two factors.
First, our method is intrinsic, so unlike the other three methods,
it does not suffer from distortions due to pulling back the mapping to a tangent space.
Second, it has clear rules for model selection,
while the other methods are often subject to model misspecification,
due to arbitrary choices of reference point, subsample points, and interpolation schemes.

\subsection{Related work}

The authors have worked on estimating functions whose domains or codomains are manifolds.
For inputs on an unknown embedded submanifold,
\cite{YangY2016} proposed a GP model that attains the minimax-optimal convergence rate,
without estimating the manifold.
To allow for noisy inputs and better scalability,
\cite{Guhaniyogi2016} first projects the input to random subspaces, and then applies a GP model.
For inputs on a known embedded submanifold,
\cite{LinLZ2019} proposed an extrinsic GP,
while \cite{NiuM2019} proposed an intrinsic GP, with heat kernel as the covariance kernel.
For outputs on an embedded submanifold,
\cite{LinLZ2017} proposed a non-GP method, which applies an extrinsic local regression
and then obtains manifold estimates via projection \cite{ZhangRD2020nr}.

While our method extends GPs to mappings that take values in the Grassmann manifold,
we are not the first to define GPs on Riemannian manifolds.
Wrapped Gaussian process (WGP) regression \cite{Mallasto2018}
approximates mappings to a general Riemannian manifold,
using wrapped Gaussian distributions (WGDs) defined by Gaussian distributions on tangent spaces.
However, this approach encounters problems when the manifold has a finite injectivity radius,
as is the case for Grassmann manifolds.
In particular, we cannot calculate the induced probability density function on the manifold
or the intrinsic mean.
In contrast, our proposed approach produces analytic forms for predictive quantities
that admit efficient computation, albeit restricted to Grassmann manifolds.

\subsection{Article structure and notations}

\Cref{sec:preliminaries} provides the basics of ROM
and the algebra and statistics of the Grassmann manifold.
\Cref{sec:subspace-regression} presents the theoretical foundation of our GPS model,
and \Cref{sec:prediction} gives an algorithm for prediction.
\Cref{sec:training-subspace} discusses model selection criteria for our model.
\Cref{sec:examples} gives several numerical experiments:
one to visualize the posterior process,
and three to access its accuracy in benchmark PROM problems.
\Cref{sec:conclusion} concludes with a discussion on practical issues.
Additional text is included in Supplementary Materials.
An \texttt{R} package accompanying this paper is available at:
\url{https://github.com/rudazhang/gpsr}.

\textit{Notations}.
Scalars are in lowercase letters, $n, k, l, d$;
vectors are in boldface lowercase letters, $\mathbf{m}, \mathbf{x}_i, \boldsymbol{\theta}$;
matrices are in boldface uppercase letters, $\mathbf{M}, \mathbf{X}_i, \mathbf{K}_l$.
Subspaces are in Fraktur script, $\mathfrak{X, M}$.
Equivalence classes are in brackets, $[\mathbf{M}], [\mathbf{m}]$.

\section{Preliminaries}
\label{sec:preliminaries}

\subsection{Parametric reduced order modeling}
\label{sub:prom}

To simplify the narrative, consider a system of ordinary differential equations (ODEs)
that is first-order, linear and time-invariant, with multiple input and output:
\begin{equation} \label{eq:system}
  \Sigma: \left\{\begin{aligned}
      \mathbf{E} \dot{\mathbf{x}} &= \mathbf{A} \mathbf{x} + \mathbf{B} \mathbf{u} \\
      \mathbf{y} &= \mathbf{C} \mathbf{x} \end{aligned}\right.
\end{equation}
With system dimension $n$, input dimension $p$, and output dimension $q$,
this system is defined by constant matrices $\mathbf{E}, \mathbf{A} \in M_{n,n}$,
$\mathbf{B} \in M_{n,p}$, and $\mathbf{C} \in M_{q,n}$.
The state $\mathbf{x}$, input $\mathbf{u}$, and output $\mathbf{y}$ are all functions of time,
with dimension $n$, $p$, and $q$ respectively.
In general, the ODE system $\Sigma$ may represent
a physical or artificial system modeled by a PDE system,
which is discretized in space, and linearized around a stationary trajectory.
The system dimension $n$ typically scales with the size of a spatial grid,
and for a large-scale problem, usually we have $n > 10^5$.

Projection-based model reduction constructs a reduced-order model (ROM) as:
\begin{equation} \label{eq:rom}
  \Sigma_r: \left\{\begin{aligned}
      \mathbf{E}_r \dot{\mathbf{x}}_r &= \mathbf{A}_r \mathbf{x}_r + \mathbf{B}_r \mathbf{u} \\
      \mathbf{y}_r &= \mathbf{C}_r \mathbf{x}_r \end{aligned}\right.
\end{equation}
Let $\mathbf{V}, \mathbf{W} \in V_{k,n}$ be orthonormal bases of $k$-dimensional subspaces,
the reduced system matrices are defined as
$\mathbf{E}_r = \mathbf{W}^T \mathbf{E} \mathbf{V}$,
$\mathbf{A}_r = \mathbf{W}^T \mathbf{A} \mathbf{V}$,
$\mathbf{B}_r = \mathbf{W}^T \mathbf{B}$, and $\mathbf{C}_r = \mathbf{C} \mathbf{V}$.
Therefore we have $\mathbf{E}_r, \mathbf{A}_r \in M_{k,k}$, $\mathbf{B}_r \in M_{k,p}$,
and $\mathbf{C}_r \in M_{q,k}$.
If the reduced bases $\mathbf{V}$ and $\mathbf{W}$ are the same,
this framework is called the Galerkin projection;
otherwise, it is called the Petrov-Galerkin projection.
Usually we would expect a reduced system dimension $k \le 50$.
Because simulation time and model storage scale at least linearly with system dimension,
they are reduced by several orders of magnitude via ROM.

To compute a reduced basis for the Galerkin projection,
a widely-used classic method is called the proper orthogonal decomposition (POD),
originally proposed for turbulent flow analysis by \cite{Lumley1967}.
This method takes a collection of system states $\mathbf{x}(t_i)$ at discrete times $\{t_i\}_{i=1}^m$,
called snapshots, which may be obtained via simulation or experimental measurements.
Let $\mathbf{X}$ be the matrix that stacks the snapshots as column vectors,
then the POD basis $\mathbf{V}$ corresponds to the left singular vectors of $\mathbf{X}$
associated with the largest $k$ singular values.
This means that the POD basis minimizes the $\mathcal{L}_2$ error of snapshot reconstruction,
which is an appealing property of POD.
For large-scale systems, the number of snapshots required is far less than the system dimension,
and usually $m = \mathcal{O}(10^3)$.

Another class of ROM methods are interpolatory \cite{Antoulas2020},
which approximate the transfer function of the original system using rational interpolation.
The transfer function of the system $\Sigma$ is defined as
$\mathbf{H}(s) = \mathbf{C} (s \mathbf{E} - \mathbf{A})^{-1} \mathbf{B}$.
Here, $\mathbf{H}: \mathbb{C} \mapsto M_{q,p}(\mathbb{C})$
is a complex matrix-valued function of a complex frequency variable.
These methods interpolate the transfer function at an arbitrary number of points
and up to an arbitrary number of derivatives along certain tangent directions.
Among such methods, the iterative rational Krylov algorithm (IRKA)
introduced by \cite{Gugercin2008} has seen great success,
which we will discuss more in \cref{sub:microthruster}.

Besides POD and interpolatory methods,
there are other frequency domain approaches such as balanced truncation \cite{Moore1981},
most common in systems and control theory,
and time domain approaches such as DMD \cite{Schmid2008},
which also discovers coherent structure in time.
There are effective ROM methods for systems more general than \cref{eq:system},
such as DEIM \cite{Chaturantabut2010} for nonlinear systems.

Our discussion so far assumes that the full model $\Sigma$ is constant.
In a more general class of problems, $\Sigma$ is parametric, such that the system matrices
$\mathbf{E}, \mathbf{A}, \mathbf{B}$, and $\mathbf{C}$ depend on a set of parameters
$\boldsymbol{\theta} \in \Theta \subset \mathbb{R}^d$.
This dependency can be nonlinear in general,
and the dimension of the parameter space is often modest, with $d \le 10$.
There are many methods for PROM, and we refer the readers to \cite{Benner2015}.
Among the four types of PROM methods discussed therein,
subspace interpolation and our proposed method belong to interpolating among local reduced bases,
while matrix and manifold interpolation belong to interpolating among local reduced system matrices.
Both types of methods require generating reduced bases at a sample of the parameter space,
computed using any projection-based ROM method.
Because generating a ROM can be computationally expensive, the sample size $l$ cannot be too large.
For problems with a modest number of parameters, usually $l \in [20, 100]$.

\subsection{Grassmann manifold}
\label{sub:Grassmann}

Because we are building a probabilistic surrogate of subspace-valued mappings,
it is helpful to review the algebra and statistics
of the Grassmann manifold and some related matrix manifolds.
For some basics of the algebra and differential geometry, see e.g. \cite{Bendokat2020,ZhangRD2021nbb};
for the statistics, see \cite{Chikuse2003}.

Let $M_{n,k}$ be the set of all $n$-by-$k$ matrices of real numbers,
which can be identified as the Euclidean space $\mathbb{R}^{n \times k}$.
The set of all full-rank $n$-by-$k$ matrices is
$M_{n,k}^* = \{\mathbf{M} \in M_{n,k} : \text{rank}(\mathbf{M}) = \min(n,k)\}$.
When $k = n$, the manifold $M_{n,k}^*$ coincides with the general linear group $\text{GL}_n$,
which consists of full-rank order-$n$ matrices.
The Stiefel manifold $V_{k,n}$ consists of all orthonormal k-frames in the Euclidean n-space:
$V_{k, n} = \{\mathbf{X} \in M_{n,k}^* : \mathbf{X}^T \mathbf{X} = \mathbf{I}_k\}$,
where $k \le n$ and $\mathbf{I}_k$ is the order-$k$ identity matrix.
The order of the subscripts is reversed by convention.
When $k = n$, the Stiefel manifold coincides with the orthogonal group $O(n)$.
Define projection $\pi: M_{n,k}^* \mapsto V_{k,n}$,
such that for any $\mathbf{M} \in M_{n,k}^*$ with a thin singular value decomposition (SVD)
$\mathbf{M} = \mathbf{V} \boldsymbol{\Sigma} \mathbf{U}^T$, $\mathbf{V} \in V_{k,n}$,
$\mathbf{U} \in O(k)$, we have $\pi(\mathbf{M}) = \mathbf{V} \mathbf{U}^T$.
Although the SVD is not unique, this mapping is uniquely defined.

The Grassmann manifold $G_{k,n}$ consists of all $k$-subspaces of the Euclidean $n$-space:
$G_{k, n} = \{\text{span}(\mathbf{M}) : \mathbf{M} \in M^*_{n, k}\}$,
where $\text{span}(\mathbf{M})$ denotes the subspace spanned by the columns of $\mathbf{M}$.
Every element of $G_{k, n}$ is a subspace, which is often represented by a basis of the subspace.
For example, every $\mathbf{M} \in M^*_{n, k}$ represents $\mathfrak{M} = \text{span}(\mathbf{M})$,
its column vectors form a basis of $\mathfrak{M}$, and every element in its equivalence class
$[\mathbf{M}] = \{\mathbf{M} \mathbf{A} : \mathbf{A} \in \text{GL}_k\}$
represents $\mathfrak{M}$ as well.
We call $\mathbf{M}$ a basis representation of $\mathfrak{M}$.
In particular, every $\mathbf{X} \in V_{k,n}$ represents $\mathfrak{X} = \text{span}(\mathbf{X})$,
and its column vectors form an orthonormal basis of $\mathfrak{X}$.
We call $\mathbf{X}$ a Stiefel representation of $\mathfrak{X}$.

The Grassmann manifold is often identified with the set of rank-$k$ symmetric projection matrices:
let $\mathcal{S}(n)$ be the set of order-$n$ symmetric matrices, define $P_{k, n} =
\{\mathbf{P} \in \mathcal{S}(n) : \mathbf{P}^2 = \mathbf{P}, \text{rank}(\mathbf{P}) = k\}$.
This identification is possible because the mapping that takes a matrix to its range
is a bijection from $P_{k, n}$ to $G_{k,n}$.
Given a Stiefel representation $\mathbf{X}$,
a subspace $\mathfrak{X}$ can be uniquely identified as $\mathbf{X} \mathbf{X}^T$.
Due to this explicit identification,
probability distributions on the Grassmann manifold can be induced through distributions on $P_{k, n}$,
with the corresponding probability density function (PDF) being:
$p: P_{k, n} \mapsto \mathbb{R}_{\ge 0}$, $\int_{P_{k, n}} p(\mathbf{P}) \mu(d \mathbf{P}) = 1$,
where $\mu$ is the normalized invariant measure on $P_{k, n}$
under the group action of $\text{GL}_{n}$.
Any probability distribution on $M_{n,k}$
that is invariant under right-orthogonal transformation
induces a probability distribution on $G_{k, n}$ \cite[Thm 2.4.8]{Chikuse2003}:
let $p$ be a PDF on $M_{n,k}$ such that $p(\mathbf{M}) = p(\mathbf{M} \mathbf{Q})$
for all $\mathbf{M} \in M_{n,k}$ and $\mathbf{Q} \in O(k)$,
if $\mathbf{M} \sim p$, let $\mathbf{X} = \pi(\mathbf{M}) \sim p_V$
and $\mathbf{X} \mathbf{X}^T \sim p_G$, then $p_V(\mathbf{X}) = p_V(\mathbf{X} \mathbf{Q})$
for all $\mathbf{Q} \in O(k)$, and $p_G(\mathbf{X} \mathbf{X}^T) = p_V(\mathbf{X})$.

Now we introduce some common probability distributions on matrix manifolds.
Let $\mathcal{S}_+(n)$ be the set of order-$n$ positive-definite matrices.
Let $\mathbf{M} \in M_{n,k}$, $\boldsymbol{\Sigma}_1 \in \mathcal{S}_+(n)$,
and $\boldsymbol{\Sigma}_2 \in \mathcal{S}_+(k)$.
The $n$-by-$k$ matrix-variate Gaussian distribution
$N_{n,k}(\mathbf{M}; \boldsymbol{\Sigma}_1, \boldsymbol{\Sigma}_2)$ is the distribution of
$\mathbf{Y} = \boldsymbol{\Sigma}_1^{1/2} \mathbf{Z} \boldsymbol{\Sigma}_2^{1/2} + \mathbf{M}$,
where $\mathbf{Z}$ is a random $n$-by-$k$ matrix
whose entries are independent standard Gaussian random variables.
The vectorized matrix $\mathbf{Y}$ is an $(nk)$-dimensional Gaussian random vector
with a special form of covariance matrix:
$\text{vec}(\mathbf{Y}) \sim N_{nk}(\text{vec}(\mathbf{M}),
\boldsymbol{\Sigma}_2 \otimes \boldsymbol{\Sigma}_1)$,
where $\text{vec}()$ denotes vectorization of a matrix by stacking its columns,
and $\otimes$ is the Kronecker product.
The matrix angular central Gaussian distribution $\text{MACG}(\boldsymbol{\Sigma})$
is a probability distribution on $V_{k,n}$, with PDF $p(\mathbf{X}; \boldsymbol{\Sigma})
= z^{-1} |\mathbf{X}^T \boldsymbol{\Sigma}^{-1} \mathbf{X}|^{-n/2}$,
where $|\cdot|$ denotes the determinant, normalizing constant $z = |\boldsymbol{\Sigma}|^{k/2}$
and parameter $\boldsymbol{\Sigma} \in \mathcal{S}_+(n)$.
This parametric family contains the uniform distribution:
since $p(\mathbf{X}; \mathbf{I}_n) = 1$, we have $\text{MACG}(\mathbf{I}_n) \sim \text{Uniform}$.
The parameter of the MACG distribution
is identified up to scaling: $\text{MACG}(\boldsymbol{\Sigma}) = \text{MACG}(c \boldsymbol{\Sigma})$,
for all $\boldsymbol{\Sigma} \in \mathcal{S}_+(n)$ and $c \in \mathbb{R}_{>0}$.
Because the MACG distribution is invariant under right-orthogonal transformation,
it also defines a family of distributions on $G_{k,n}$ with the same PDF.
These distributions are related to the matrix-variate Gaussian as follows:
let $\mathbf{M} \sim N_{n,k}(0; \boldsymbol{\Sigma}, \mathbf{I}_k)$
where $\boldsymbol{\Sigma} \in \mathcal{S}_+(n)$;
let $\mathbf{X} = \pi(M)$, then $\mathbf{X} \sim \text{MACG}(\boldsymbol{\Sigma})$,
and $\mathbf{X} \mathbf{X}^T \sim \text{MACG}(\boldsymbol{\Sigma})$.
Due to the above property,
$\text{MACG}(\boldsymbol{\Sigma})$ can be easily sampled
by generating $\mathbf{M} \sim N_{n,k}(0; \boldsymbol{\Sigma}, \mathbf{I}_k)$
and projecting it via $\pi$.

\section{Gaussian process subspace regression}
\label{sec:subspace-regression}

We now present the proposed Gaussian Process Subspace regression (GPS) model.
Because GPs take values in a Euclidean space,
they are not directly applicable to approximate a subspace-valued mapping $f: \Theta \mapsto G_{k,n}$,
where the codomain is the Grassmann manifold.
Instead, we may find vector-valued mappings $\bar{f}: \Theta \mapsto \mathbb{R}^{nk}$
that are representations of $f$,
in the sense that $f = \text{span} \circ \text{vec}^{-1} \circ \bar{f}$.
Here, $\circ$ denotes the composition of two mappings
and $\text{vec}^{-1}: \mathbb{R}^{nk} \mapsto M_{n,k}$
denotes the ``inverse'' of $\text{vec}()$, that is,
constructing a matrix columnwise from a vector.
Such representations are not unique, and we denote the set of representations as
$\bar{F} = \{\bar{f} : f = \text{span} \circ \text{vec}^{-1} \circ \bar{f}\}$.
Now $f$ can be identified with $\bar{F}$,
or equivalently, any distribution supported on $\bar{F}$.

GP models extend naturally to approximate distributions on a set of functions.
Let $\mathfrak{X} = f(\boldsymbol{\theta})$ with a basis representation $\mathbf{X}$.
Recall that $\mathbf{X}$ has an equivalence class
$[\mathbf{X}] = \{\mathbf{X} \mathbf{A} : \mathbf{A} \in \text{GL}_k\}$.
Let $\mathbf{x} = \text{vec}(\mathbf{X})$.
Its equivalence class can be written as
$[\mathbf{x}] = \{\text{vec}(\mathbf{X} \mathbf{A}) : \mathbf{A} \in \text{GL}_k\}$.
Assume that $\bar{f}$ have a GP prior, we may assign equal likelihood to $[\mathbf{x}]$.
We can then proceed to derive the posterior and the predictive distributions.
In the following, we provide modeling details and analytical solutions for this approach.

\subsection{Model specification}
\label{sub:model}

We start by specifying a prior for the representations.
Without other information on $f$, an uninformative prior is
for $f(\boldsymbol{\theta})$ to be uniformly distributed on $G_{k,n}$.
We can achieve this by assigning $\bar{f}(\boldsymbol{\theta}) \sim N_{nk}(0, \mathbf{I}_{nk})$,
the $nk$-dimensional standard Gaussian.
To see this, let matrix $\mathbf{M} = \text{vec}^{-1}(\bar{f}(\boldsymbol{\theta}))$,
then $\mathbf{M} \sim N_{n,k}(0; \mathbf{I}_n, \mathbf{I}_k)$ is a matrix-variate standard Gaussian;
let subspace $\mathfrak{M} = \text{span}(\mathbf{M})$,
then $\mathfrak{M} \sim \text{MACG}(\mathbf{I}_n) \sim \text{Uniform}(G_{k,n})$.
We assign a correlation structure as follows.
Let $k: \Theta \times \Theta \mapsto [-1, 1]$ be a correlation function,
i.e. a positive definite kernel with $k(\boldsymbol{\theta}, \boldsymbol{\theta}) = 1$
for all $\boldsymbol{\theta} \in \Theta$.
For any finite collection of input points $\boldsymbol{\theta} = (\boldsymbol{\theta}_i)_{i=1}^l$,
let $\mathbf{m}_i = \bar{f}(\boldsymbol{\theta}_i)$,
and let $\mathbf{K}_l$ be the order-$l$ correlation matrix with entry
$[\mathbf{K}_l]_{ij} = k(\boldsymbol{\theta}_i, \boldsymbol{\theta}_j)$.
We assign the function values $\mathbf{m} = (\mathbf{m}_i)_{i=1}^l$ a prior joint distribution
$\mathbf{m} \sim N_{nkl}(0, \mathbf{K}_l \otimes \mathbf{I}_{nk})$.
Compactly, we can write this GP prior as $\bar{f} \sim \mathcal{GP}(0, k \otimes \mathbf{I}_{nk})$.
This is the simplest covariance structure for $\bar{f}$.

Without a likelihood function, this GP prior gives predictions as follows.
Let $\boldsymbol{\theta}_*$ be a target point and $\mathbf{m}_* = \bar{f}(\boldsymbol{\theta}_*)$.
We have the prior joint distribution:
\begin{equation} \label{eq:prior-subspace}
  (\mathbf{m}_*, \mathbf{m}) \sim N_{nk(l+1)}(0, \mathbf{K}_{l+1} \otimes \mathbf{I}_{nk})
\end{equation}
where $\mathbf{K}_{l+1} = [1 \; \mathbf{k}_l^T; \mathbf{k}_l \; \mathbf{K}_l]$
and $\mathbf{k}_l = (k(\boldsymbol{\theta}_*, \boldsymbol{\theta}_i))_{i=1}^l$.
If we write $\mathbf{K}_{22} = \mathbf{K}_l \otimes \mathbf{I}_{nk}$
and $\mathbf{K}_{12} = \mathbf{k}_l^T \otimes \mathbf{I}_{nk}$,
by properties of multivariate Gaussian distributions,
the conditional distribution of $\mathbf{m}_*$ given $\mathbf{m}$ can be written as:
\begin{align} \label{eq:predictive-point-subspace}
  \mathbf{m}_* | \mathbf{m}
  &\sim N_{nk}(\mathbf{K}_{12} \mathbf{K}_{22}^{-1} \mathbf{m},
    \mathbf{I}_{nk} - \mathbf{K}_{12} \mathbf{K}_{22}^{-1} \mathbf{K}_{12}^T) \\
  &= N_{nk}\left(\sum_{i=1}^l [\mathbf{K}_l^{-1} \mathbf{k}_l]_i \mathbf{m}_i,
    (1- \mathbf{k}_l^T \mathbf{K}_l^{-1} \mathbf{k}_l) \mathbf{I}_{nk}\right) \nonumber
\end{align}

We assign equal likelihood to the equivalence class of representations.
Assume that we have function evaluations $\mathfrak{X}_i = f(\boldsymbol{\theta}_i)$
with Stiefel representations $\mathbf{X}_i \in V_{k,n}$.
Let $\mathbf{x}_i = \text{vec}(\mathbf{X}_i)$
and $[\mathbf{x}_i] = \{\text{vec}(\mathbf{X}_i \mathbf{A}) : \mathbf{A} \in \text{GL}_k\}$.
For $\mathbf{m}_i = \bar{f}(\boldsymbol{\theta}_i)$, the likelihood function gives:
\begin{equation} \label{eq:likelihood-subspace}
  L(\mathbf{m}_i | \mathfrak{X}_i) = 1(\mathbf{m}_i \in [\mathbf{x}_i])
\end{equation}

The posterior distribution of $\mathbf{m}$ given observations
$\mathfrak{X} = (\mathfrak{X}_i)_{i=1}^l$ is derived from the prior and the likelihood via Bayes’ rule:
\begin{equation} \label{eq:posterior-subspace}
  p(\mathbf{m} | \mathfrak{X}) \propto
  \exp\bigg\{-\frac{1}{2} \mathbf{m}^T (\mathbf{K}_l \otimes \mathbf{I}_{nk})^{-1} \mathbf{m}\bigg\}
  \prod_{i=1}^l 1(\mathbf{m}_i \in [\mathbf{x}_i])
\end{equation}

\subsection{Predictive distributions}
\label{sub:predictive}

The predictive distribution of $\mathbf{m}_*$ given observations $\mathfrak{X}$
is obtained by integrating the conditional distribution \cref{eq:predictive-point-subspace}
over the posterior distribution \cref{eq:posterior-subspace}.
We summarize the result as follows:
\begin{theorem} \label{thm:predictive-subspace}
Let $\mathbf{X} = [\mathbf{X}_1~\cdots~\mathbf{X}_l]$
be the matrix that combines $\mathbf{X}_i$ by columns,
and $\mathbb{X} = \diag(\mathbf{X}_i)_{i=1}^l$
be the matrix with $\mathbf{X}_i$ as diagonal blocks.
Let $\varepsilon^2 = 1 - \mathbf{k}_l^T \mathbf{K}_l^{-1} \mathbf{k}_l$,
$\mathbf{v} = \mathbf{K}_l^{-1} \mathbf{k}_l$, $\mathbf{D}_{\mathbf{v}} = \diag(\mathbf{v})$,
and $\widetilde{\mathbf{K}}_l = (\mathbf{D}_{\mathbf{v}} \mathbf{K}_l \mathbf{D}_{\mathbf{v}})^{-1}$.
The predictive distribution of $\mathbf{m}_*$ given observations $\mathfrak{X}$ is:
\begin{gather}
  \mathbf{m}_* | \mathfrak{X}
  \sim N_{nk}(0, \mathbf{I}_k \otimes \boldsymbol{\Sigma}) \nonumber \\
  \boldsymbol{\Sigma} = \varepsilon^2 \mathbf{I}_n + \mathbf{X}
  [\mathbb{X}^T (\widetilde{\mathbf{K}}_l \otimes \mathbf{I}_n) \mathbb{X}]^{-1} \mathbf{X}^T
  \label{eq:covariance-subspace}
\end{gather}
\end{theorem}
The proof is quite lengthy and thus deferred to \cref{app:proof}.
This theorem shows that, given observations:
(1) the matrix $\mathbf{M}_* = \text{vec}^{-1}(\mathbf{m}_*)$
has a matrix-variate Gaussian distribution
$\mathbf{M}_* | \mathfrak{X} \sim N_{n,k}(0; \boldsymbol{\Sigma}, \mathbf{I}_k)$;
and (2) the subspace $\mathfrak{M}_* = \text{span}(\mathbf{M}_*)$
has an MACG distribution $\mathfrak{M}_* | \mathfrak{X} \sim \text{MACG}(\boldsymbol{\Sigma})$
(see \cref{sub:Grassmann}).

The predictive distributions admit an intuitive interpretation.
Since $\boldsymbol{\Sigma}$ is positive semi-definite,
there is an eigenvalue decomposition (EVD)
$\boldsymbol{\Sigma} = \mathbf{Q} \diag(\boldsymbol{\lambda}) \mathbf{Q}^T$,
where $\boldsymbol{\lambda} \in \mathbb{R}_{\ge 0}^n$ are in decreasing order
and $\mathbf{Q} \in O(n)$.
Therefore we can simulate $\mathbf{M}_* | \mathfrak{X}$ as
$\mathbf{M}_* = \boldsymbol{\Sigma}^{1/2} \mathbf{Z} =
\mathbf{Q} \diag(\boldsymbol{\lambda})^{1/2} \mathbf{Q}^T \mathbf{Z}$,
where $\mathbf{Z} \in M_{n,k}$ is a random matrix of standard Gaussians.
The column vectors of $\mathbf{Z}$ are scaled by the square root of the eigenvalue in each eigenspace;
therefore the range of $\mathbf{M}_*$ is more likely to align with
the top eigenspaces of $\boldsymbol{\Sigma}$.
Recall that $\mathfrak{M}_* = \text{span}(\mathbf{M}_*)$. We have the following results.
The global Riemannian center of mass of $\mathfrak{M}_* | \mathfrak{X}$
is $\text{span}(\mathbf{V})$, where $\mathbf{V}$ is the first $k$ columns of $\mathbf{Q}$.
The uncertainty of $\mathfrak{M}_* | \mathfrak{X}$
is compactly described by the eigenvalues $\boldsymbol{\lambda}$:
the larger an eigenvalue is, the more important is the associated eigenspace;
and the mean prediction is more useful if $(\lambda_i)_{i=k+1}^n$ are small
relative to $(\lambda_i)_{i=1}^k$.

A main feature of our GP model is that,
while its construction involves the extrinsic Euclidean space of basis representations of subspaces,
its predictive distribution is intrinsic to the Grassmann manifold.
In particular, our model does not involve tangent spaces or the Riemannian exponential,
and thus it is not subject to the distortions associated with applying local tangent approximations.
Moreover, the function space explored by the GPS is much broader
than existing interpolation approaches,
so our model is more flexible and robust to model misspecification.
Perhaps surprisingly, the GPS has closed-form expressions for its predictive distributions,
which enables efficient computation for subspace prediction and uncertainty quantification.

While \cref{thm:predictive-subspace} is concerned with point predictions on the Grassmann manifold,
our GP model also induces joint distributions on $G_{k,n}$
and can be used to generate random subspace-valued functions
(see \cref{sup:joint-distribution}).

\section{Prediction algorithm}
\label{sec:prediction}

From \Cref{thm:predictive-subspace} and the discussion thereafter we see that,
to compute the predictive distribution, one needs the EVD of $\boldsymbol{\Sigma}$.
Here we give an efficient method to compute this.

\subsection{Efficient EVD of $\boldsymbol{\Sigma}$}

Denote $\boldsymbol{\Pi} = \mathbb{X}^T (\widetilde{\mathbf{K}}_l \otimes \mathbf{I}_n) \mathbb{X}$
and $\check{\boldsymbol{\Sigma}} = \mathbf{X} \boldsymbol{\Pi}^{-1} \mathbf{X}^T$.
We note that $\widetilde{\mathbf{K}}_l$, $\boldsymbol{\Pi}$, and
$\check{\boldsymbol{\Sigma}}$ are all positive semi-definite.
Let $r = \text{rank}(\mathbf{X}) \le \min(n, k l)$,
then $\check{\boldsymbol{\Sigma}}$ also has rank $r$ and therefore $r$ positive eigenvalues.
From the form of $\check{\boldsymbol{\Sigma}}$,
we see that its top-$r$ eigenvectors span the range of $\mathbf{X}$.
Let $\mathbf{X} = \widetilde{\mathbf{V}} \widetilde{\mathbf{R}} \widetilde{\mathbf{P}}^T$
be a rank-revealing QR decomposition,
such that $\widetilde{\mathbf{V}}$ has $r$ orthonormal columns,
$\widetilde{\mathbf{R}}$ is an $r$-by-$kl$ upper triangular matrix,
and $\widetilde{\mathbf{P}}$ is a permutation matrix.
Denote $\mathbf{S} = \widetilde{\mathbf{V}}^T \check{\boldsymbol{\Sigma}} \widetilde{\mathbf{V}}$
and let $\mathbf{S} =
\mathring{\mathbf{Q}} \diag(\mathring{\boldsymbol{\lambda}}) \mathring{\mathbf{Q}}^T$
be an EVD where $\mathring{\boldsymbol{\lambda}}$ is descending and
$\mathring{\mathbf{Q}} \in O(r)$.
Let $\mathbf{V} = \widetilde{\mathbf{V}} \mathring{\mathbf{Q}}$
and let $\mathbf{Q} = (\mathbf{V}, \mathbf{V}_\perp) \in O(n)$ be an orthogonal completion.
Let $\check{\boldsymbol{\lambda}} = (\mathring{\boldsymbol{\lambda}}, \mathbf{0}_{n-r})$
where $\mathbf{0}_{n-r}$ is the vector of zeros with length $n-r$.
Then we have an EVD:
$\check{\boldsymbol{\Sigma}} = \mathbf{Q} \diag(\check{\boldsymbol{\lambda}}) \mathbf{Q}^T$.  
Because $\boldsymbol{\Sigma} = \check{\boldsymbol{\Sigma}} + \varepsilon^2 \mathbf{I}_n$,
we have an EVD of $\boldsymbol{\Sigma}$:
\begin{equation} \label{eq:covariance-EVD}
  \boldsymbol{\Sigma} = \mathbf{Q} \diag(\check{\boldsymbol{\lambda}} +
  \varepsilon^2 \mathbf{1}_n) \mathbf{Q}^T
\end{equation}
Here $\mathbf{1}_n$ is the vector of ones with length $n$.
We see that, for a complete probabilistic prediction, we only need
a rank-revealing QR of $\mathbf{X}$, an EVD of $\mathbf{S}$, and $\varepsilon^2$.
For the mean prediction, we only need the top-$k$ eigenvectors of $\mathbf{S}$.

We can simplify the computation of $\mathbf{S}$ as follows.
Note that $\widetilde{\mathbf{V}}^T \mathbf{X} = \widetilde{\mathbf{R}} \widetilde{\mathbf{P}}^T$
and $\widetilde{\mathbf{P}}^{-1} = \widetilde{\mathbf{P}}^T$.
Because $\mathbf{S} = \widetilde{\mathbf{V}}^T \check{\boldsymbol{\Sigma}} \widetilde{\mathbf{V}}$
and $\check{\boldsymbol{\Sigma}} = \mathbf{X} \boldsymbol{\Pi}^{-1} \mathbf{X}^T$,
we have $\mathbf{S} = \widetilde{\mathbf{R}}
(\widetilde{\mathbf{P}} \boldsymbol{\Pi} \widetilde{\mathbf{P}}^T)^{-1} \widetilde{\mathbf{R}}^T$.
Let Gram matrix $\Box = \mathbf{X}^T \mathbf{X}$, which has a block matrix structure
$\Box = [\Box_{ij}]_{i,j=1}^l$ with $\Box_{ij} = \mathbf{X}_i \mathbf{X}_j$.
Note that $\boldsymbol{\Pi}$ similarly has a block matrix structure
$\boldsymbol{\Pi} = [\boldsymbol{\Pi}_{ij}]_{i,j=1}^l$
with $\boldsymbol{\Pi}_{ij} = \widetilde{k}_{ij} \Box_{ij}$,
where $\widetilde{k}_{ij} = [\widetilde{\mathbf{K}}_l]_{i,j}$.
The construction of $\boldsymbol{\Pi}$ can be written in a compact form:
$\boldsymbol{\Pi} = \Box \circ (\widetilde{\mathbf{K}}_l \otimes \mathbf{J}_k)$,
where $\circ$ denotes the Hadamard product and $\mathbf{J}_k$ is the order-$k$ matrix of ones.
Let $\widetilde{\boldsymbol{\Pi}} = \widetilde{\mathbf{P}} \boldsymbol{\Pi} \widetilde{\mathbf{P}}^T$
and let $\widetilde{\boldsymbol{\Pi}} = \mathbf{L} \mathbf{L}^T$ be a Cholesky decomposition,
where $\mathbf{L}$ is a lower triangular matrix.
Let $\widetilde{\mathbf{L}} = \mathbf{L}^{-1} \widetilde{\mathbf{R}}^T$ by solving linear equations,
which is a $kl$-by-$r$ lower triangular matrix,
then we have $\mathbf{S} = \widetilde{\mathbf{L}}^T \widetilde{\mathbf{L}}$.

We formally describe the prediction procedure in two parts:
\Cref{alg:GPS-preprocess} only needs to be done once,
and \Cref{alg:GPS-predict-EVD} is needed for each prediction.

\alglanguage{pseudocode}
\begin{algorithm}[h]
  \caption{GPS: Preprocessing}
  \label{alg:GPS-preprocess}
  \begin{algorithmic}[1] %
    \Input observation $\mathbf{X} = [\mathbf{X}_1~\cdots~\mathbf{X}_l]$.
    \State Compute Gram matrix: $\Box \gets \mathbf{X}^T \mathbf{X}$.
    \label{line:Gram}
    \State Rank-revealing QR:
    $\mathbf{X} = \widetilde{\mathbf{V}} \widetilde{\mathbf{R}} \widetilde{\mathbf{P}}^T$.
    \label{line:QR}
    \Output Gram matrix $\Box$; global basis $\widetilde{\mathbf{V}}$;
    upper triangular $\widetilde{\mathbf{R}}$; pivoting $\widetilde{\mathbf{P}}$.
  \end{algorithmic}
\end{algorithm}

\alglanguage{pseudocode}
\begin{algorithm}[h]
  \caption{GPS: Prediction}
  \label{alg:GPS-predict-EVD}
  \begin{algorithmic}[1] %
    \Require correlation function $k(\cdot, \cdot)$;
    preprocessing output
    $(\Box, \widetilde{\mathbf{V}}, \widetilde{\mathbf{R}}, \widetilde{\mathbf{P}})$.
    \Input sample $(\boldsymbol{\theta}_i)_{i=1}^l$;
    target $\boldsymbol{\theta}_*$;
    truncation size $t \in \{k, k+1, \cdots, r\}$.
    \State Construct correlation matrix and vector:
    $k_{ij} \gets k(\boldsymbol{\theta}_i, \boldsymbol{\theta}_j)$,
    $k_i \gets k(\boldsymbol{\theta}_*, \boldsymbol{\theta}_i)$.
    \label{line:K}
    \State Solve linear equations: $\mathbf{v} \gets \text{solve}(\mathbf{K}, \mathbf{k})$,
    $\widehat{\mathbf{K}} \gets \text{solve}(\mathbf{K}, \diag(\mathbf{v})^{-1})$.
    \label{line:small-solve}
    \State Construct matrix: $\boldsymbol{\Pi} \gets [\boldsymbol{\Pi}_{ij}]_{i,j=1}^l$, where
    $\boldsymbol{\Pi}_{ij} \gets v_i^{-1} \widehat{k}_{ij} \Box_{ij}$.
    \label{line:dummy-prod}    
    \State Cholesky decomposition:
    $\widetilde{\mathbf{P}} \boldsymbol{\Pi} \widetilde{\mathbf{P}}^T = \mathbf{L} \mathbf{L}^T$.
    \label{line:Cholesky}
    \State Solve linear equations:
    $\widetilde{\mathbf{L}} \gets \text{solve}(\mathbf{L}, \widetilde{\mathbf{R}}^T)$ 
    \label{line:solve-triangular}
    \State Cross product: $\mathbf{S} \gets \widetilde{\mathbf{L}}^T \widetilde{\mathbf{L}}$.
    \label{line:cross-product}
    \State Truncated EVD: $\mathbf{S} = \mathring{\mathbf{V}}
    \diag(\mathring{\boldsymbol{\lambda}}) \mathring{\mathbf{V}}^T$,
    where $\mathring{\boldsymbol{\lambda}}$ has length $t$.
    \label{line:truncated-EVD}
    \State Compute noise variance: $\varepsilon^2 \gets 1 - \mathbf{k}^T \mathbf{v}$.
    \label{line:noise-var}
    \Output principal directions $\mathbf{V} = \widetilde{\mathbf{V}} \mathring{\mathbf{V}}$;
    principal variances $\mathring{\boldsymbol{\lambda}}$; noise variance $\varepsilon^2$.
    \Note May return $\widetilde{\mathbf{V}}$ and $\mathring{\mathbf{V}}$
    instead of $\mathbf{V}$ to avoid matrix multiplication.
  \end{algorithmic}
\end{algorithm}

\subsection{Computational cost and comparison}

Here we analyze the computational cost of each step in floating point operations (flops),
accurate up to the dominant term.
In \Cref{alg:GPS-preprocess},
line~\ref{line:Gram} takes $n k^2 l^2$ flops;
line~\ref{line:QR} takes $\mathcal{O}(n kl r)$ flops, and if $r \approx kl$, this requires
about $4 n k^2 l^2$ flops using the Householder QR with column pivoting \cite{Golub2013}.
In \Cref{alg:GPS-predict-EVD},
line~\ref{line:K} evaluates the correlation function $l^2 / 2$ times;
line~\ref{line:small-solve} takes $l^3 / 3$ flops for Cholesky decomposition,
and $2 l^3$ for forward and back substitution;
line~\ref{line:dummy-prod} takes $k^2 l^2 / 2$ flops;
line~\ref{line:Cholesky} takes $k^3 l^3 / 3$ flops;
line~\ref{line:solve-triangular} takes $k^3 l^3 / 3 - (kl - r)^3 / 3$ flops,
due to the upper triangular structure in $\widetilde{\mathbf{R}}$;
line~\ref{line:cross-product} takes $r^3 / 3 + (kl - r) r^2$ flops,
due to the lower triangular structure in $\widetilde{\mathbf{L}}$;
line~\ref{line:truncated-EVD} takes $\mathcal{O}(r^2 t)$
with classical or randomized algorithms \cite{Halko2011};
and line~\ref{line:noise-var} takes $2 l$ flops.
Note that $\mathbf{K}$ and its Cholesky decomposition can be reused for future predictions.
Overall, with $n > kl$ and assuming $r \approx kl$ and $t = k$,
\Cref{alg:GPS-preprocess} gives an overhead cost of
about $5 n k^2 l^2$ flops if we use the Householder QR with column pivoting,
and \Cref{alg:GPS-predict-EVD} gives a cost of about $k^3 l^3$ flops per prediction.

An alternative version of \Cref{alg:GPS-predict-EVD}
is to conduct a truncated singular value decomposition (SVD): $\widetilde{\mathbf{L}}
= \mathring{\mathbf{V}} \diag(\mathring{\boldsymbol{\sigma}}) \mathbf{W}^T$,
and then return $\mathring{\mathbf{V}}$ and
$\mathring{\boldsymbol{\lambda}} = \mathring{\boldsymbol{\sigma}}^2$.
Although this avoids the cross product in line~\ref{line:cross-product}
and thus saves about $k^3 l^3 / 3$ flops,
truncated SVD can take a significant amount of time and eliminate the saving.
Theoretically, the truncated SVD takes $\mathcal{O}(r kl t)$ with classical algorithms,
and $\mathcal{O}(r kl \log t)$ with randomized algorithms \cite{Halko2011}.
But in practice, the truncated SVD appears to be more costly than the truncated EVD.
Since truncated SVD gives a less accurate result than truncated EVD,
we consider \Cref{alg:GPS-predict-EVD} as the reference version.

Note that the matrix multiplication $\mathbf{V} = \widetilde{\mathbf{V}} \mathring{\mathbf{V}}$
takes $2 n r t$ flops, which would dominate the cost per prediction if $n > k l^2 / 2$.
However, this cost can be avoided if principal directions $\mathbf{V}$ are not explicitly needed.
In PROM problems, to compute an order-$k$ reduced matrix
$\mathbf{A}_k = \mathbf{V}^T \mathbf{A} \mathbf{V}$,
one may precompute an order-$r$ matrix
$\mathbf{A}_r = \widetilde{\mathbf{V}}^T \mathbf{A} \widetilde{\mathbf{V}}$,
and then compute $\mathbf{A}_k = \mathring{\mathbf{V}}^T \mathbf{A}_r \mathring{\mathbf{V}}$.
Since $\mathbf{A}$ is usually sparse,
the cost of a matrix-vector multiplication $\mathbf{A} \mathbf{x}$
is usually $T_{\text{mult}} = \mathcal{O}(n)$.
Then the computation has an overhead cost of $2 n k^2 l^2 + kl T_{\text{mult}}$ flops,
and only takes about $2 k^3 l^2$ flops per prediction.

In comparison, subspace interpolation \cite{Amsallem2008}
takes $\mathcal{O}(n k^2 n_r)$ flops per prediction,
where $n_r \le l$ is the total number of sample points used
and $\mathcal{O}(n k^2)$ corresponds to the cost of a thin SVD and matrix multiplications.
Usually $n_r$ is set to a small number, about 4 or 5.
Since its prediction does not have a special factorization structure,
to compute a reduced matrix it takes another $2 n k^2 + k T_{\text{mult}}$ flops,
where $T_{\text{mult}}$ denotes the cost of a matrix-vector multiplication.
The prediction cost can be greatly reduced if the problem has only one parameter
and one uses linear interpolation \cite{Son2013}.
As mentioned in \cref{sub:literature},
matrix interpolation \cite{Panzer2010} and manifold interpolation \cite{Amsallem2011}
directly interpolate local ROMs so their prediction costs do not depend on $n$,
and therefore they are suitable for online computation.
\Cref{tab:flops} compares the computational costs of the proposed method
and these three interpolatory methods for PROM problems.

Our method is typically much faster than methods for computing local reduced bases.
Consider the computation of a local POD basis given $m$ snapshots at one parameter point.
The cost is dominated by a truncated SVD of the n-by-m snapshot matrix,
which takes $\mathcal{O}(n m k)$ time.
To compare the costs, take the rocket injector example in \cite{Mak2018},
where $n \approx 10^5$, $m = 10^3$, $k = 45$, $l = 30$.
We have $(n m k) / (k^3 l^3) \approx 1.83$.
Considering the constant factor in truncated SVD,
in this case our method is about an order of magnitude faster than computing a local POD basis.
Because the cost of computing snapshots dominates the overall POD procedure,
this implies a clear advantage in using our method to approximate local POD bases.
The cost of computing a pair of local IRKA bases is less straightforward to analyze
\cite{Antoulas2020,Hokanson2020}.
Every iteration needs to solve $2k$ systems of linear equations,
each with a different coefficient matrix of order $n$ that cannot be reused across iterations.
The number of iterations depends on the initial values provided to the algorithm,
and the algorithm needs to be restarted
if it does not converge after a predefined maximum number of iterations.
Depending on the problem, IRKA can take longer than the POD procedure.

\begin{table*}[tb]
  \centering
  \begin{threeparttable}
    \label{tab:flops}
    \caption{Interpolatory methods for PROM: flop counts of the dominant terms.}
    \begin{tabular}{lccccp{0.2\linewidth}}
      \toprule
      & Preprocess & Subspace & ROM & Tuning & Reference \\
      \midrule
      GPS & $5 n k^2 l^2$ & $k^3 l^3$ & $2 k^3 l^2$ & $k^3 l^4$ & this paper \\
      Subspace-Int & $10 n k^2 l^2$ & $8 n k^2$ & $2 n k^2$ & \textdagger & \cite{Amsallem2008} \\
      Matrix-Int & $6 n k^2 l^2$ & - & $2 k^2 l$ & \textdagger & \cite{Panzer2010}   \\
      Manifold-Int & $\mathcal{O}(n k^2 l)$ & - & $\mathcal{O}(k^3 l)$\textsuperscript{*}
                                    & \textdagger & \cite{Amsallem2011} \\
      \bottomrule
    \end{tabular}
    \begin{tablenotes}
    \item[*] Coefficient usually on the scale of 50 due to matrix exponential / logarithm,
      which can be numerically unstable \cite{Higham2005}.
    \item[\textdagger] Optimal choice of reference ROM and interpolation scheme
      is an open problem.
    \end{tablenotes}
  \end{threeparttable}
\end{table*}

\section{Model selection}
\label{sec:training-subspace}

Although the correlation function $k(\cdot, \cdot)$ can be arbitrary,
it is often specified in a form that depends on some hyperparameters \cite[Ch. 4]{Rasmussen2006}.
For example, the squared exponential (SE) kernel is:
\begin{equation} \label{eq:SE-kernel}
  k(\boldsymbol{\theta}, \boldsymbol{\theta}'; \boldsymbol{\beta}) =
  \prod_{i=1}^d \exp[- \frac{1}{2} (\theta_i-\theta_i')^2 / \beta_i^2]
\end{equation}
where length-scales $\boldsymbol{\beta} = (\beta_i)_{i=1}^d$ are the hyperparameters.
GP models with the SE kernel are smooth, and the length-scales can be understood as
characteristic distances along each parameter before the function values become uncorrelated.

One can set the hyperparameters to optimize a certain criterion to improve prediction,
see e.g. \cite[Sec 5.4]{Rasmussen2006} and \cite[Sec 3.3]{Santner2018}.
For GPS, we recommend minimizing the leave-one-out cross validation (LOOCV) predictive error,
measured in Riemannian distances.
(Other distances between subspaces may be used as well,
but we choose Riemannian distance for concreteness.)
In this section we analyze and give an algorithm to compute this criterion.
We provide a procedure to compute its gradient in \cref{sup:gradient}
and discuss some alternative criteria in \cref{sup:criteria}.

In our experience, the predictive performance of GPS
is not very sensitive to hyperparameters,
so one may use certain default values to trade accuracy for reduced computational cost.
As a rule-of-thumb for the SE kernel,
one may set the length-scales to $3 d^{3/2}/ l$ relative to the parameter ranges,
and expect good predictive results.

\subsection{LOOCV predictive error}

To measure predictive error, we need a score of dissimilarity for pairs of subspaces.
There are many metrics defined on the Grassmann manifold, see e.g. \cite{YeK2016} for a list.
Among them, the most commonly used is the Riemannian distance,
which is the length of the shortest curves connecting two points in a Riemannian manifold.
The Riemannian distance between subspaces $\mathfrak{X}, \mathfrak{Y} \in G_{k,n}$
is the 2-norm of their principal angles, which can be computed as \cite{Bendokat2020}:
\begin{equation}\label{eq:Riemannian-distance}
  d_g(\mathfrak{X}, \mathfrak{Y}) = \|\arccos \boldsymbol{\sigma}(\mathbf{X}^T \mathbf{Y})\|
\end{equation}
Here, $\mathbf{X}, \mathbf{Y} \in V_{k,n}$ are representations of the subspaces,
and $\boldsymbol{\sigma}(\cdot)$ denotes the singular values of a matrix.
Let $\mathbf{V}_{-i}$ represent the mean prediction for target $\boldsymbol{\theta}_i$,
using observations $(\boldsymbol{\theta}_j, \mathbf{X}_j)_{j \ne i}$.
The LOOCV predictive error can be defined as:
\begin{equation}\label{eq:LOOCV-error}
  \epsilon_2 = \sum_{i=1}^l d_g^2(\mathbf{X}_i, \mathbf{V}_{-i}) =
  \sum_{i=1}^l \sum_{j=1}^k \left(\arccos \sigma_j(\mathbf{X}_i^T \mathbf{V}_{-i})\right)^2
\end{equation}
Here we use the sum of squared errors for its smoothness and,
with a slight abuse of notation, we replace the subspaces with their representations.

We can compute $\mathbf{X}_i^T \mathbf{V}_{-i}$ efficiently.
Note that $\mathbf{V}_{-i}$ consists of the top-$k$ eigenvectors of $\boldsymbol{\Sigma}_{-i}$
which, analogous to \cref{eq:covariance-subspace}, can be written as:
\begin{equation} \label{eq:covariance-subspace-LOO}
  \boldsymbol{\Sigma}_{-i} = \varepsilon_{-i}^2 \mathbf{I}_n + \mathbf{X}_{-i}
  [\mathbb{X}_{-i}^T (\widetilde{\mathbf{K}}_{-i} \otimes \mathbf{I}_n) \mathbb{X}_{-i}]^{-1}
  \mathbf{X}_{-i}^T
\end{equation}
Here, all the quantities are defined without the $i$-th observation.
Similar to the analysis in \cref{sec:prediction},
denote $\boldsymbol{\Pi}_{-i} =
\mathbb{X}_{-i}^T (\widetilde{\mathbf{K}}_{-i} \otimes \mathbf{I}_n) \mathbb{X}_{-i}$ and
$\check{\boldsymbol{\Sigma}}_{-i} = \mathbf{X}_{-i} (\boldsymbol{\Pi}_{-i})^{-1} \mathbf{X}_{-i}^T$.
Let $r_{-i} = \text{rank}(\mathbf{X}_{-i})$, then the top-$r_{-i}$ eigenvectors of
$\check{\boldsymbol{\Sigma}}_{-i}$ span the range of $\mathbf{X}_{-i}$,
which is a subset of the range of $\mathbf{X}$.
Recall that $\mathbf{X} = \widetilde{\mathbf{V}} \widetilde{\mathbf{R}} \widetilde{\mathbf{P}}^T$
is a rank-revealing QR, let
$\mathbf{S}_{-i} = \widetilde{\mathbf{V}}^T \check{\boldsymbol{\Sigma}}_{-i} \widetilde{\mathbf{V}}$
and let $\mathring{\mathbf{V}}_{-i}$ be the top-$k$ eigenvectors of $\mathbf{S}_{-i}$,
then $\widetilde{\mathbf{V}} \mathring{\mathbf{V}}_{-i}$
are the top-$k$ eigenvectors of $\check{\boldsymbol{\Sigma}}_{-i}$,
which are the same as those of $\boldsymbol{\Sigma}_{-i}$.
Hence, $\mathbf{V}_{-i} = \widetilde{\mathbf{V}} \mathring{\mathbf{V}}_{-i}$.
Let $\widetilde{\mathbf{C}} = \widetilde{\mathbf{R}} \widetilde{\mathbf{P}}^T
= \widetilde{\mathbf{V}}^T \mathbf{X}$
and let $\widetilde{\mathbf{C}}_{i} = \widetilde{\mathbf{V}}^T \mathbf{X}_{i}$,
which can be obtained simply by subsetting $\widetilde{\mathbf{C}}$
with column indices from $(i-1)k+1$ to $ik$.
Then we have $\mathbf{X}_i^T \mathbf{V}_{-i} =
\mathbf{X}_i^T \widetilde{\mathbf{V}} \mathring{\mathbf{V}}_{-i} =
\widetilde{\mathbf{C}}_i^T \mathring{\mathbf{V}}_{-i}$.
Similarly, let $\widetilde{\mathbf{C}}_{-i} = \widetilde{\mathbf{V}}^T \mathbf{X}_{-i}$,
which can be obtained by removing $\widetilde{\mathbf{C}}_{i}$ from $\widetilde{\mathbf{C}}$,
then we have $\mathbf{S}_{-i} = \widetilde{\mathbf{C}}_{-i}
(\boldsymbol{\Pi}_{-i})^{-1} \widetilde{\mathbf{C}}_{-i}^T$.

To compute $\boldsymbol{\Pi}_{-i}$, note that $\widetilde{\mathbf{K}}_{-i} =
(\mathbf{D}_{\mathbf{v}_{-i}} \mathbf{K}_{-i} \mathbf{D}_{\mathbf{v}_{-i}})^{-1}$,
where $\mathbf{D}_{\mathbf{v}_{-i}} = \diag(\mathbf{v}_{-i})$,
$\mathbf{v}_{-i} = (\mathbf{K}_{-i})^{-1} \mathbf{k}_{-i}$,
$\mathbf{K}_{-i} = [k_{pq}]_{p,q \ne i}$, and $\mathbf{k}_{-i} = (k_{pi})_{p \ne i}$.
The LOOCV quantities $\mathbf{v}_{-i}$ and $(\mathbf{K}_{-i})^{-1}$
can be written in terms of $\mathbf{K}^{-1}$,
see for example \cite[Sec. 5.4.2]{Rasmussen2006}.
Let $\bar{\mathbf{K}} = \mathbf{K}^{-1}$,
then $\mathbf{v}_{-i} = -\bar{k}_{ii}^{-1} \bar{\mathbf{k}}_{-i}$ and
$(\mathbf{K}_{-i})^{-1} = \bar{\mathbf{K}}_{-i} - \bar{k}_{ii} \mathbf{v}_{-i} \mathbf{v}_{-i}^T$.
Here, $\bar{\mathbf{K}}_{-i} = [\bar{k}_{pq}]_{p,q \ne i}$
and $\bar{\mathbf{k}}_{-i} = (\bar{k}_{pi})_{p \ne i}$.
Together we have $\widetilde{\mathbf{K}}_{-i} = \bar{k}_{ii}^{-1} \Delta_{-i}$,
where $\Delta_{-i} = [\bar{k}_{pq} \bar{k}_{ii} / (\bar{k}_{ip} \bar{k}_{iq}) - 1]_{p,q \ne i}$.
Also, $\boldsymbol{\Pi}_{-i}$ can be written in a compact form:
$\boldsymbol{\Pi}_{-i} = \bar{k}_{ii}^{-1} \Box_{-i} \circ (\Delta_{-i} \otimes \mathbf{J}_k)$.
Since we are only concerned with the eigenvectors of $\mathbf{S}_{-i}$,
with a little abuse of notation,
we redefine $\boldsymbol{\Pi}_{-i}$ without the term $\bar{k}_{ii}^{-1}$.
We describe the overall procedure in \Cref{alg:GPS-LOOCV-error}.

\alglanguage{pseudocode}
\begin{algorithm}[h]
  \caption{LOOCV Predictive Error}
  \label{alg:GPS-LOOCV-error}
  \begin{algorithmic}[1] %
    \Require correlation function $k$; %
    sample $(\boldsymbol{\theta}_i)_{i=1}^l$;
    preprocessing output
    $(\Box, \widetilde{\mathbf{C}} = \widetilde{\mathbf{R}} \widetilde{\mathbf{P}}^T)$.    
    \Input hyperparameters $\boldsymbol{\beta}$.
    \State Construct inverse correlation matrix:
    $\bar{\mathbf{K}} \gets \text{solve}(\mathbf{K})$,
    where $k_{ij} \gets k(\boldsymbol{\theta}_i, \boldsymbol{\theta}_j; \boldsymbol{\beta})$.
    \label{line:K-inv}
    \For{$i$ in $1, \cdots, l$}
    \State Construct: $\boldsymbol{\Pi} \gets [\boldsymbol{\Pi}_{pq}]_{p,q \ne i}$, where
    $\boldsymbol{\Pi}_{pq} \gets \delta_{pq} \Box_{pq}$,
    $\delta_{pq} \gets \bar{k}_{pq} \bar{k}_{ii} / (\bar{k}_{ip} \bar{k}_{iq}) - 1$.
    \label{line:LOO-Pi}
    \State Construct:
    $\mathbf{S} \gets \widetilde{\mathbf{L}}^T \widetilde{\mathbf{L}}$,
    where $\boldsymbol{\Pi} = \mathbf{L} \mathbf{L}^T$,
    $\widetilde{\mathbf{L}} \gets \text{solve}(\mathbf{L}, \widetilde{\mathbf{C}}_{-i}^T)$.
    \label{line:LOO-S}    
    \State Truncated EVD: $\mathbf{S} = \mathring{\mathbf{V}}
    \diag(\mathring{\boldsymbol{\lambda}}) \mathring{\mathbf{V}}^T$,
    where $\mathring{\boldsymbol{\lambda}}$ has length $k$.
    \label{line:LOO-truncated-EVD}
    \State Compute singular values:
    $\boldsymbol{\sigma} \gets \boldsymbol{\sigma}(\widetilde{\mathbf{C}}_i^T \mathring{\mathbf{V}})$.
    \label{line:singular-values}
    \State Compute squared error: $\epsilon_i \gets \sum_{j=1}^k \arccos(\sigma_j)^2$
    \label{line:error}
    \EndFor
    \Output LOOCV predictive error $\epsilon = \sum_{i=1}^l \epsilon_i$.
  \end{algorithmic}
\end{algorithm}

\subsection{Computational cost}

In terms of computational cost, \Cref{alg:GPS-LOOCV-error}
is approximately $l$ repetitions of \Cref{alg:GPS-predict-EVD},
so it costs about $k^3 l^4$ flops per evaluation.
Since evaluating the LOOCV error and gradient (see \cref{sup:gradient})
is about $l$ times the prediction cost,
hyperparameter tuning may be a significant part of the overall cost.
In practice, we recommend setting a very rough convergence threshold;
for parameters with a range of one, a threshold of 0.01 is sufficient for the length-scale.
If the problem has multiple parameters, once they are scaled into comparable ranges,
they may share the same length-scale.
If multiple hyperparameters are to be trained,
gradient-based optimization methods can be more efficient than just using the LOOCV error.
To minimize the number of iterations, one may also set a restrictive range and, if applicable,
a good initial value for the hyperparameters;
for example, $\pm 30\%$ of the aforementioned rule-of-thumb length-scale,
with initial value at the midpoint.

\section{Numerical experiments}
\label{sec:examples}

\subsection{Visualization of GP subspace regression}
\label{sub:viz}

\begin{figure}[t]
  \centering
  \includegraphics[width=0.8\linewidth]{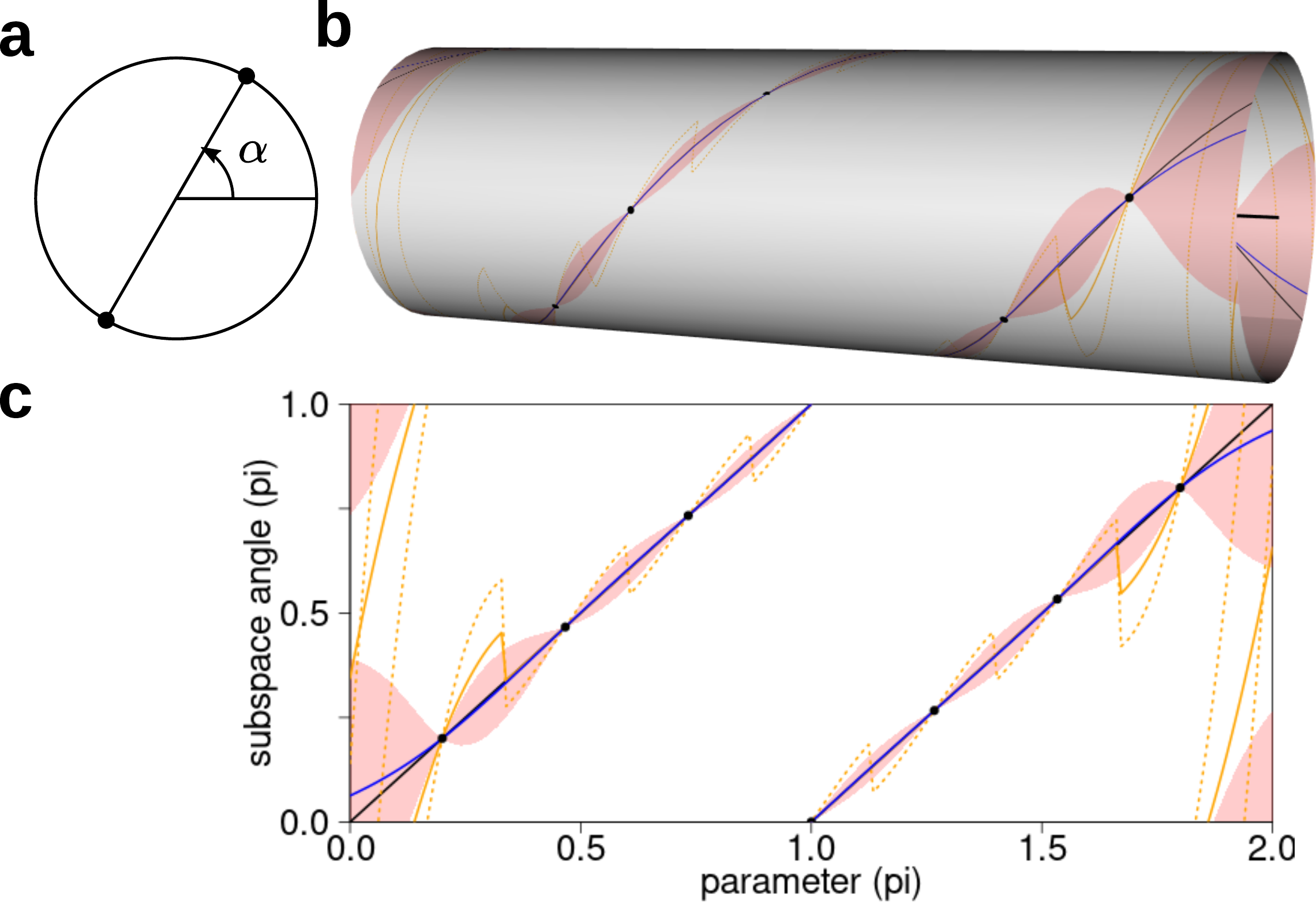}
  \caption{Visualization of the GPS model.
    (\textbf{a}) Every 1d subspace in the plane can be uniquely identified
    by either a pair of antipodal points on a circle, or an angle $\alpha \in [0, \pi)$.
    (\textbf{b}) Posterior process of the GPS model on the surface of a cylinder:
    true function (black line), data (black points),
    predictive mean (blue curve), 95\% predictive interval (red shade).
    Orange curves are predictions from subspace interpolation: $n_r = 3$ (solid), $n_r = 4$ (dotted).
    (\textbf{c}) Same as (\textbf{b}) but as a 2d plot.}
  \label{fig:visualization}
\end{figure}

The simplest type of subspace-valued functions have the form $f: \mathbb{R} \mapsto G_{1,2}$,
which maps a real number to a one-dimensional linear subspace in the plane.
The Grassmann manifold $G_{1,2}$ can be identified as the unit circle,
treating antipodal points as equivalent (\Cref{fig:visualization}a).
Therefore, such a function $f$ can be plotted on the surface of a cylinder
(\Cref{fig:visualization}b),
which helps us visualize the posterior process of GP subspace regression.

Specifically, let $f$ be a covering map such that $f(\theta)$
is the subspace with angle $\alpha = \theta~\text{mod}~\pi$.
This can be plotted as a double helix on the cylinder.
To approximate this function with the proposed GPS model,
suppose we observe sample points $\theta_i = c_i \pi$,
where $c_i$ are seven equal-distanced points between 0.2 and 1.8.
For the correlation function $k$, we use the SE kernel,
and set the length-scale $\beta$ by minimizing the LOOCV predictive error.
In this example, $\beta = 2.8 \approx 0.9 \pi$.
To visualize predictive uncertainty,
we plot the 95\% posterior predictive intervals from Theorem \ref{thm:predictive-subspace}.
We also include results from the existing subspace interpolation method \cite{Amsallem2008}
for comparison.
As suggested by the authors, for every target parameter
we use the nearest $n_r$ sampled points for the interpolation
(where $n_r = 3$ and $4$ in \Cref{fig:visualization}),
among which the nearest sampled point is used as the reference point.
We use Lagrange interpolation for the tangent vectors.

We see that, with only seven data points, the predictive mean function of GPS
closely tracks the true function within the range of sampled parameter points.
Furthermore, the uncertainties from our model also well-cover the truth:
the posterior predictive intervals contain the true subspace values for all $\theta \in [0, 2 \pi]$.
Note that as the target point moves away from the sample points,
the predictive distribution degenerates to the prior, the uniform distribution on $G_{1,2}$.
Subspace interpolation, on the other hand, yields noticeably poorer predictions
compared to GPS for both $n_r=3$ and $n_r = 4$.
As a deterministic interpolation approach,
it also does not provide a quantification of interpolation uncertainty.
This shows that, for this example, the proposed GPS model uses sample data more effectively
to yield better predictions with uncertainty quantification.

\subsection{Anemometer: approximating local POD bases}
\label{sub:anemometer}

Here we consider a benchmark problem for PROM known as the anemometer \cite{anemometer},
a type of micro-electromechanical system (MEMS) device
that measures the flow speed of its surroundings.
Such a device needs to be calibrated under different flow conditions for its temperature response.
However, an accurate representation of the device needs to resolve
the coupled fluid and thermodynamics, and can be very time-consuming to compute.
It is therefore useful to apply PROM methods.

Specifically, a convection-diffusion equation
is discretized into a linear ODE system as \cref{eq:system},
with system dimension $n = 29,008$ and input and output dimensions $p = q = 1$.
The matrix $\mathbf{A}$ depends on one parameter $\theta \in [0, 1]$ representing fluid velocity
and is not symmetric in general,
while $\mathbf{E}, \mathbf{B}, \mathbf{C}$ are constants.
The input map $\mathbf{B}$ represents a heat source,
and the output map $\mathbf{C}$ gives the temperature difference of two nodes.

To build a parametric reduced-order model (PROM),
we first construct local POD bases at a sample of the parameter space,
and then use the mean prediction of GPS
to estimate the reduced subspaces at other parameter points.
As before, we use the SE kernel, with a length-scale
that minimizes the LOOCV predictive error.
The subspace-valued mappings being approximated in this problem
have very high dimensional codomains:
because the dimension of $G_{k,n}$ is $k (n - k)$,
with $k=20$ and $k=40$, the manifold dimensions here are 579,760 and 1,158,720 respectively.

For comparison, we also estimate the reduced subspaces using subspace interpolation,
with the same setup as in the visualization example.
For manifold interpolation \cite{Amsallem2011},
we use the same setup for subspace interpolation.
For matrix interpolation \cite{Panzer2010},
we use the nearest sampled point as the reference point
and, as suggested by the authors,
we use linear interpolation for the reduced system matrices.
We include results for local POD bases as a reference level we would like to match.

To measure the error introduced by a ROM, a standard choice is the $\mathcal{H}_2$ metric,
defined as the largest possible amplitude of the output error given any unit-energy input.
Let $\|\cdot\|_{\mathcal{L}_2}$ and $\|\cdot\|_{\mathcal{L}_\infty}$ denote
the $\mathcal{L}_2$ and $\mathcal{L}_\infty$ norms, respectively.
Following the notations of \cref{eq:system,eq:rom}, we have:
\begin{equation} \label{eq:H2-metric}
  \|\Sigma - \Sigma_r\|_{\mathcal{H}_2} =
  \sup_{\mathbf{u} \in \mathcal{L}_2} \frac{\|\mathbf{y} - \mathbf{y}_r\|_{\mathcal{L}_\infty}}
  {\|\mathbf{u}\|_{\mathcal{L}_2}}
\end{equation}
Relative $\mathcal{H}_2$ error is the $\mathcal{H}_2$ error
divided by the $\mathcal{H}_2$ norm of the original system:
\begin{equation} \label{eq:relative-H2}
  e(\Sigma, \Sigma_r)_{\mathcal{H}_2}
  = \frac{\|\Sigma - \Sigma_r\|_{\mathcal{H}_2}}{\|\Sigma\|_{\mathcal{H}_2}}
\end{equation}
The $\mathcal{H}_2$ norms can be obtained analytically via the controllability Gramian,
which can be computed by solving the Lyapunov equations \cite{mmess}.

\begin{figure}[t]
  \centering
  \includegraphics[width=0.98\linewidth]{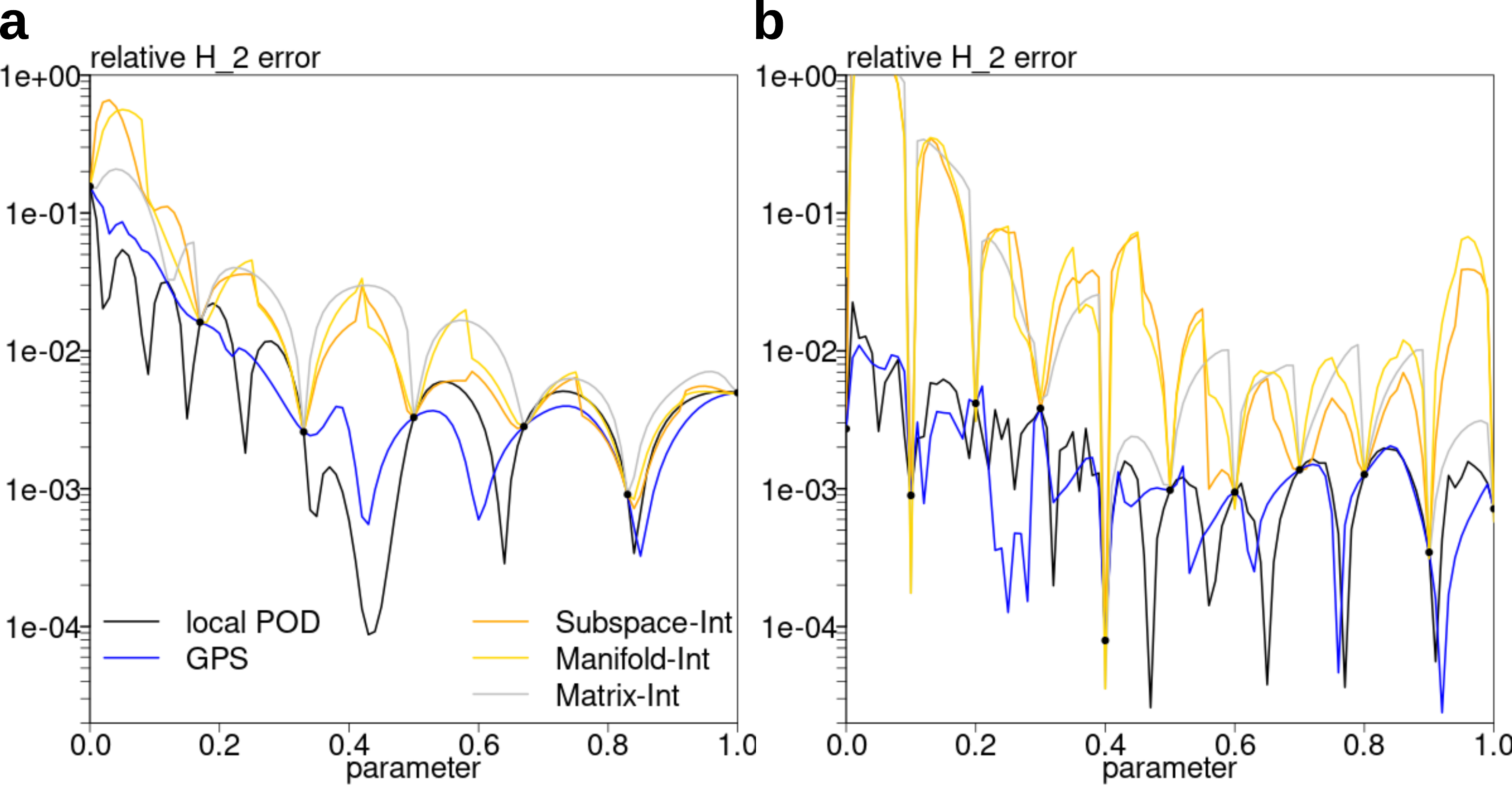}
  \caption{Anemometer, relative $\mathcal{H}_2$ error:
    (\textbf{a}) $k=20$; (\textbf{b}) $k=40$.
    Training data shown as points.
  }
  \label{fig:anemometer}
\end{figure}

\Cref{fig:anemometer}a shows the relative $\mathcal{H}_2$ errors using these methods,
with subspace dimension $k = 20$.
Here we use a sample of seven equal-distanced points from 0 to 1.
GPS uses a length-scale $\beta = 0.36$, selected via LOOCV.
The results for subspace and manifold interpolation use $n_r = 3$;
the results are similar for $n_r = 4$ or $5$.
We see that the three existing interpolation methods perform similarly,
and the errors tend to blow up in between sample points.
In comparison, the proposed GPS model yields much lower errors:
the relative $\mathcal{H}_2$ error is comparable to that for the local POD (the reference level).
Note that the goal here is not to perfectly match the error curve of local POD,
but to keep the error as low as possible;
in this sense, the GPS model appears to provide noticeable improvements over existing methods.

\Cref{fig:anemometer}b shows the results for $k = 40$.
Here we use a sample of 11 equal-distanced points from 0 to 1.
GPS uses a length-scale $\beta = 0.25$.
Setup for the interpolation methods are unchanged.
We see that, even with the increased sample size,
all three interpolation methods fail to keep a low error level.
While matrix interpolation occasionally does better than the other two,
this is probably not generalizable due to the linear interpolation scheme.
In comparison, our method again yields much lower errors,
and maintains a similar level of accuracy as the local POD.

Another error measure is the $\mathcal{L}_2$ state error.
The $\mathcal{L}_2$ metric of square-integrable functions on the interval $[0, T]$,
discretized into $J$ parts of length $\delta t$, can be approximated as:
\begin{equation} \label{eq:L2-metric}
  \|\mathbf{x} - \hat{\mathbf{x}}\|_{\mathcal{L}_2}^2
  = \int_0^T \|\mathbf{x}(t) - \hat{\mathbf{x}}(t)\|_2^2~dt
  \approx \sum_{i=1}^J \|\mathbf{x}(t_i) - \hat{\mathbf{x}}(t_i)\|_2^2~\delta t
\end{equation}
Relative $\mathcal{L}_2$ state error is the $\mathcal{L}_2$ error of the state vector of a ROM,
divided by the $\mathcal{L}_2$ norm of the state vector of the original system.
Following %
\cref{eq:system,eq:rom}, this gives:
\begin{equation} \label{eq:relative-L2}
  e(\mathbf{x}, \mathbf{x}_r)_{\mathcal{L}_2}
  = \frac{\|\mathbf{x} - \mathbf{V} \mathbf{x}_r\|_{\mathcal{L}_2}}{\|\mathbf{x}\|_{\mathcal{L}_2}}
\end{equation}

\begin{figure}[t]
  \centering
  \includegraphics[width=0.98\linewidth]{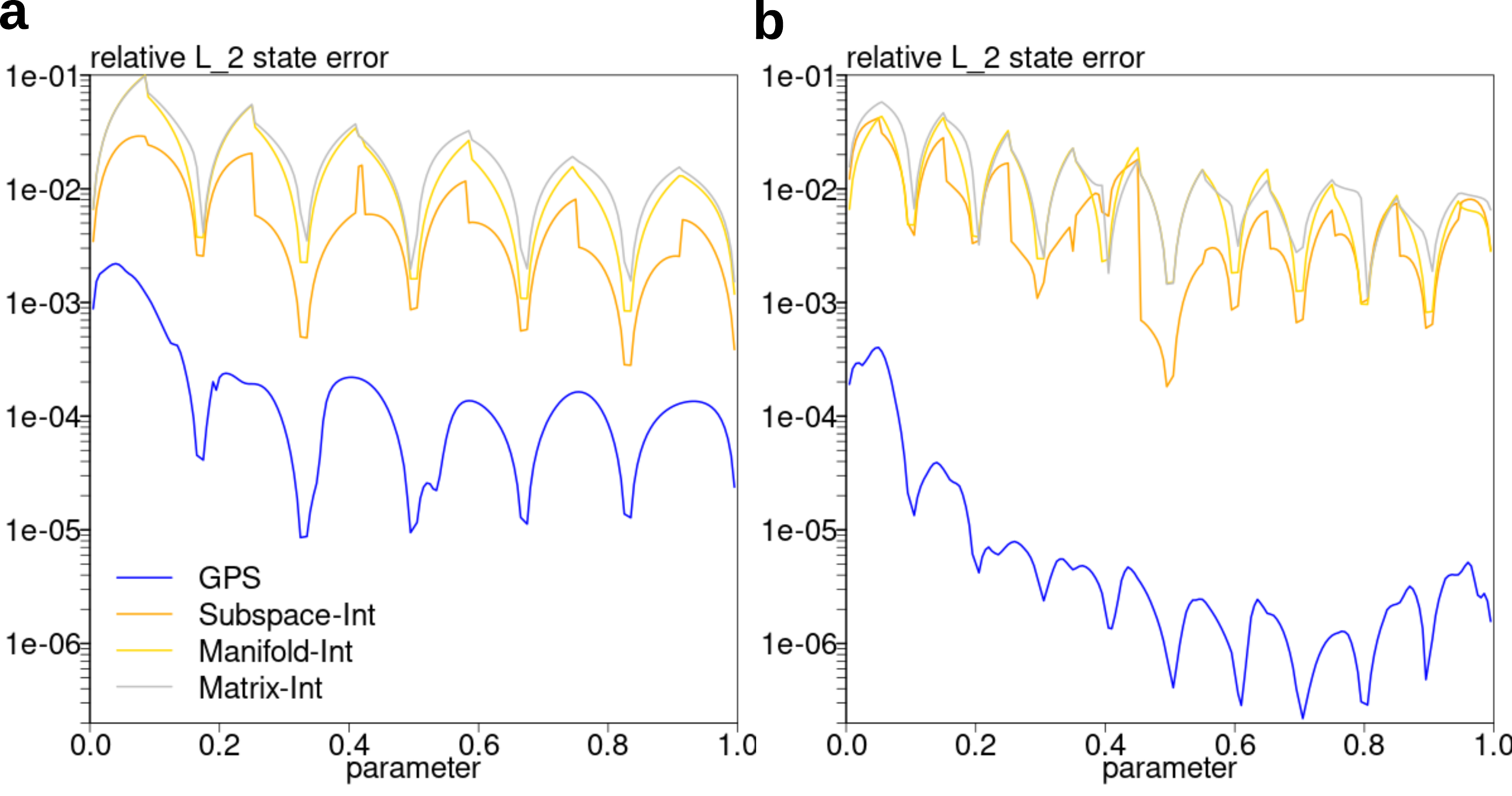}
  \caption{Anemometer, relative $\mathcal{L}_2$ state error:
    (\textbf{a}) $k=20$; (\textbf{b}) $k=40$.
  }
  \label{fig:anemometer-L2}
\end{figure}

\Cref{fig:anemometer-L2} shows the relative $\mathcal{L}_2$ state errors using these methods.
Local POD is omitted from these plots
since its relative $\mathcal{L}_2$ state error is practically zero.
The error curves of the three interpolation methods are qualitatively similar,
with subspace interpolation better than manifold interpolation,
which is in turn better than matrix interpolation.
In comparison, the GPS again yields much lower errors:
for $k = 20$, the average error is about two orders of magnitude lower
than that of subspace interpolation;
for $k = 40$, it is about three orders of magnitude lower.
This improvement can be attributed to the more flexible and intrinsic nature of the GPS model,
which allows for more effective use of sample data.

\subsection{Microthruster: approximating local IRKA bases}
\label{sub:microthruster}

Here we consider another benchmark problem for PROM known as the microthruster \cite{thermalmodel},
an array of solid propellant microthrusters on a chip.
To find an optimal design of array geometry and driving circuit,
many simulations need to be carried out, which can be prohibitive with large-scale models.
The use of PROM is therefore justified.

Specifically, the numerical model discretizes a heat transfer equation
into a linear ODE system as \cref{eq:system},
with system dimension $n = 4,257$, input dimension $p = 1$, and output dimension $q = 7$.
The input $\mathbf{B}$ represents the electrical circuit,
and the output $\mathbf{C}$ gives the temperature at seven nodes.
The convection boundary conditions are parameterized into three parameters,
each within the range $[1, 10^4]$,
and affect the symmetric system matrix $\mathbf{A}$ on the diagonal.
To simplify comparison, we fix the three parameters to always be the same,
and take the base-10 logarithm of their original values,
so we have one parameter $\theta \in [0, 4]$.

For this problem, we use IRKA to construct reduced bases at the sample points.
IRKA \cite{Gugercin2008,Antoulas2020} is an iterative algorithm
that searches for an order-$k$ rational function that approximates the transfer function,
until it satisfies the tangential interpolation conditions.
If IRKA converges, the converged point locally minimizes the $\mathcal{H}_2$ error
in the space of order-$k$ rational functions.
IRKA constructs a ROM in state space via the two-sided Petrov-Galerkin projection,
that is, the reduced bases $\mathbf{V}$ and $\mathbf{W}$ are different.

Because IRKA uses two different bases $\mathbf{V}$ and $\mathbf{W}$,
for a parametric system this means that each parameter is associated with a pair of subspaces,
and we may construct a PROM by approximating a mapping for the form
$(\mathfrak{V}, \mathfrak{W})(\boldsymbol{\theta})$.
Since our proposed method only handles mappings that output one subspace,
we proceed by modeling the pair of subspaces separately.
This inevitably leaves some information in the data unused,
and there may be methods that can improve upon this work-around.
Setup for the interpolation methods are the same as in the anemometer example.

\begin{figure}[t]
  \centering
  \includegraphics[width=0.9\linewidth]{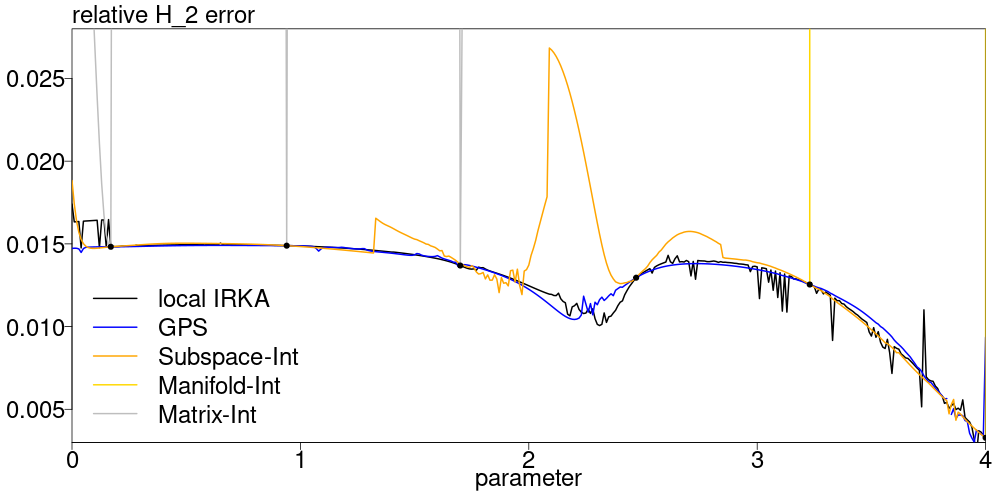}
  \caption{Microthruster, relative $\mathcal{H}_2$ error. $k=10$.
    Training data are shown as points.
  }
  \label{fig:microthruster}
\end{figure}

\Cref{fig:microthruster} shows the relative $\mathcal{H}_2$ errors using these methods,
with subspace dimension $k=10$.
Here we use a sample of 6 points: $\theta = 0.17, 0.94, 1.7, 2.47, 3.23, 4$.
GPS uses a length-scale $\beta = 1.4$ for basis $\mathbf{V}$,
and $\beta = 2.56$ for basis $\mathbf{W}$.
The result for subspace interpolation uses $n_r = 3$;
the other values of $n_r$ give results with larger errors.
We see that, while subspace interpolation matches the error curve of local IRKA
(the reference level) quite well in some parts of the parameter space,
its error blows up in an unsmooth region in between.
These errors are noticeably larger for manifold and matrix interpolation,
so we cropped them out of the plot.
To contrast, the proposed GPS method instead tracks the local IRKA error curve smoothly
across the parameter space, yielding much lower errors than existing interpolation methods.

For this problem, many of the ROMs generated by manifold interpolation are complex-valued,
due to the matrix logarithm that computes the tangent vectors.
Moreover, many ROMs generated by manifold and matrix interpolation are unstable,
which means that the $\mathcal{H}_2$ errors are infinite.
Although our method and subspace interpolation do not guarantee the stability of reduced models,
because they seem to accurately approximate the reduced subspaces,
unstable ROMs appear less often.
We discuss issues specific to approximating IRKA bases in \cref{sup:irka}.

\subsection{Anemometer: 3-parameter case}
\label{sub:anemometer-3p}

To compare the methods in a PROM problem with multiple parameters,
here we consider the three-parameter version of the anemometer \cite{anemometer}.
The parameters include specific heat $c \in [0, 1]$,
thermal conductivity $\kappa \in [1, 2]$, and fluid velocity $v \in [0.1, 2]$.
The system matrices have the form $\mathbf{E} = \mathbf{E}_s + c \mathbf{E}_f$
and $\mathbf{A} = \mathbf{A}_{d,s} + \kappa \mathbf{A}_{d,f} + c v \mathbf{A}_c$,
while $\mathbf{B}$ and $\mathbf{C}$ are constant.
Other aspects of the problem are unchanged.

To sample the parameter space,
we first use the maximin Latin hypercube sampling (LHS) to obtain a training set,
and then use the sequential maximin design to obtain a testing set,
see e.g. \cite[Ch. 4]{Gramacy2020}.
Maximin LHS generates a random set of points
that are spread out in the parameter space and well-distanced from each other.
Sequential maximin design generates another with similar properties,
but also well-distanced from the given training set.

The setup for the PROM methods remain unchanged from the 1-parameter case,
except the interpolation scheme for the three interpolation methods.
Since Lagrange and linear interpolations do not apply to multiple parameters,
we use the radial basis function (RBF) method described in \cite[p. 278]{Amsallem2010}.
Specifically, a multiquadric RBF is applied entrywise to interpolate
the tangent vectors in subspace and manifold interpolation
as well as the matrices in matrix interpolation.
For subspace interpolation,
horizontal projection is applied to maintain validity of the interpolated tangent vector.

\Cref{tab:h2-error} compares the mean relative $\mathcal{H}_2$-errors,
with training sample sizes $l = 14, 18$, or $21$, and testing sample size 100.
For each training sample, GPS uses a length-scale $\beta = 1.05, 0.85$, or $0.7$, respectively.
Notice that, with subspace dimension $k = 20$,
the mean relative $\mathcal{H}_2$-error of local POD is about 5.5\%, which is not particularly low.
It is clear that our method is able to maintain the error level of local POD
with as few as 18 training points.
In comparison, the error increase in subspace interpolation is several times higher in all cases.
Manifold interpolation is much less accurate than the previous two methods,
while matrix interpolation is the least accurate.

Similarly, \Cref{tab:l2-error} compares the mean relative $\mathcal{L}_2$ state errors.
Manifold and matrix interpolation are excluded because they cannot reconstruct the state vector.
Since local POD minimizes the $\mathcal{L}_2$ state error by construction,
its error level is practically zero.
With $l = 18$, our method has a relative error of about 1\%, half that of subspace interpolation.
This ratio drops as sample size gets larger.
Overall, the GPS method is much more data efficient than subspace interpolation
in this multi-parameter setting, again owing to its flexibility and intrinsic nature.

\begin{table*}[tb]
  \centering
  \begin{threeparttable}
    \label{tab:h2-error}
    \caption{Mean relative $\mathcal{H}_2$-error for 3-parameter anemometer,
      $k = 20$, varying sample size.}
    \begin{tabular}{lllllll}
      \toprule
      & $l = 14$ & & $l = 18$ & & $l = 21$ & \\
      \midrule
      local POD & 5.55\% & (1) & 5.46\% & (1) & 5.69\% & (1) \\
      GPS & 6.49\% & (1.169) & 5.80\% & (1.062) & 5.14\% & (0.903) \\
      Subspace-Int & 8.34\% & (1.503) & 7.38\% & (1.352) & 6.19\% & (1.148) \\
      Manifold-Int & 16.6\% & (2.986) & 13.8\% & (2.524) & 12.7\% & (2.232) \\
      Matrix-Int & 49.7\% & (8.962) & 44.2\% & (8.104) & 45.5\% & (8.003) \\
      \bottomrule
    \end{tabular}
    \begin{tablenotes}
    \item[*] Number in parentheses shows value relative to local POD.
    \end{tablenotes}
  \end{threeparttable}
\end{table*}

\begin{table*}[tb]
  \centering
  \begin{threeparttable}
    \label{tab:l2-error}
    \caption{Mean relative $\mathcal{L}_2$ state error for 3-parameter anemometer,
      $k = 20$, varying sample size.}
    \begin{tabular}{lllllll}
      \toprule
      & $l = 14$ & & $l = 18$ & & $l = 21$ & \\
      \midrule
      local POD & 7.98e-13 & (0) &  8.36e-13 & (0) & 8.77e-13 & (0) \\
      GPS & 1.24e-2 & (0.437) & 6.42e-3 & (0.273) & 5.55e-3 & (0.250) \\
      Subspace-Int & 2.85e-2 & (1) & 2.35e-2 & (1) & 2.22e-2 & (1)    \\
      \bottomrule
    \end{tabular}
    \begin{tablenotes}
    \item[*] Number in parentheses shows value relative to subspace interpolation.
    \end{tablenotes}
  \end{threeparttable}
\end{table*}

\section{Concluding remarks}
\label{sec:conclusion}

In this paper we propose a new GP model for probabilistic approximation of subspace-valued functions.
A key application of this model is for parametric reduced order modeling.
We show that the GPS model gives accurate predictions even with small sample sizes,
and because its prediction cost does not depend on system dimension $n$,
it is typically faster than subspace interpolation in PROM problems.

There are some intuitive explanations for the success of our method.
Interpolation on tangent spaces of Riemannian manifolds, such as subspace and manifold interpolation,
works best when a local interpolation scheme is performed on a few points.
When using high-degree polynomials,
these methods are subject to the oscillation between interpolation points,
a problem known as Runge's phenomenon.
Moreover, as explained in \cref{sup:discussion},
when points further away from the reference point are used,
the true mapping becomes more distorted on the tangent space and harder to approximate.
A similar concern is addressed in \cite{Zimmermann2020} Sec. 3.
Therefore, these methods cannot use more than a handful of points at a time,
and have limited potential to extend to higher-dimensional parameter spaces.
Our method is intrinsic to the Grassmann manifold and does not suffer from such limitations,
so its accuracy improves with sample size.
But since GP models are data efficient, a small sample size can still give accurate results.
As for computational cost, because subspace interpolation uses
the Riemannian exponential and logarithm of the Grassmann manifold,
both involving a thin SVD of an $n$-by-$k$ matrix,
it is slow for large-scale problems.
Our method turns the truncated EVD of an order-$n$ matrix into one of an order-$kl$ matrix,
and the prediction cost is instead dominated by construction of the matrix,
which is carried out efficiently via matrix decomposition and linear solvers.
Although manifold and matrix interpolation also scale independently of $n$,
their accuracy can degrade quickly which makes the results unusable.

Since the prediction cost of our method is cubic in subspace dimension $k$ and sample size $l$,
it is best to keep both of them small to allow for fast computation.
To keep $k$ small, one needs to choose a ROM method
that is best suited for the relevant error measure.
For example, POD is optimal in $\mathcal{L}_2$ norm of snapshot reconstruction error,
while IRKA is locally optimal in $\mathcal{H}_2$ norm of model reduction error.
To keep $l$ small, one needs to choose an efficient method for parameter sampling.
One may consider adaptive sampling and sparse grid methods \cite{Benner2015},
or experimental design methods in statistics \cite{Santner2018,Gramacy2020}.
If $l$ has to be really large, there are some GP methods that cap the $l^3$ scaling.
One approach is to use local approximate GP \cite{Gramacy2015},
where for each target point only a subsample of mostly nearby points are used in the prediction.
Another approach is covariance tapering or compactly supported kernels \cite{Furrer2006,Kaufman2011},
where the kernel becomes zero beyond a certain distance, so that the covariance matrix is sparse
and sparse matrix algorithms can be used to speed up computation.

Besides computational efficiency, another important issue in PROM
is the preservation of system properties, such as stability, passivity, and contractivity.
From the numerical examples we see that,
although stability is not guaranteed for the reduced models generated by our method,
it is still observed in most cases,
simply because our method can accurately approximate the subspace map of local ROMs.

\appendix

\section{Proof of \cref{thm:predictive-subspace}}
\label{app:proof}

We see that the posterior $p(\mathbf{m} | \mathfrak{X})$
in \cref{eq:posterior-subspace} takes positive values in $\prod_{i=1}^l [\mathbf{x}_i]$,
where $[\mathbf{x}_i] = \{\text{vec}(\mathbf{X}_i \mathbf{A}) : \mathbf{A} \in \text{GL}_k\}$.
Because $\text{GL}_k$ is a full-measure subset of $M_{k,k}$, we can replace $[\mathbf{x}_i]$
with $\{\text{vec}(\mathbf{X}_i \mathbf{A}) : \mathbf{A} \in M_{k,k}\}$
without changing the posterior.
Note that the latter equals
$\mathfrak{X}_i^k = \prod_{j=1}^k \{\mathbf{X}_i \mathbf{c} : \mathbf{c} \in \mathbb{R}^k\}$,
so the support of the posterior can be written as: $S = \prod_{i=1}^l \mathfrak{X}_i^k$.

The predictive distribution of $\mathbf{m}_*$ given observations $\mathfrak{X}$
is obtained by integrating the conditional distribution \cref{eq:predictive-point-subspace}
over the posterior distribution \cref{eq:posterior-subspace}, that is:
\begin{equation*} \label{eq:predictive-likelihood}
  \circledast := p(\mathbf{m}_* | \mathfrak{X}) =
  \int_S p(\mathbf{m}_* | \mathbf{m})~p(\mathbf{m} | \mathfrak{X})~d \mathbf{m}
\end{equation*}
Every $\mathbf{m} \in S$ can be written as $\mathbf{m} = (\mathbf{m}_i)_{i=1}^l$,
where $\mathbf{m}_i = \text{vec}(\mathbf{X}_i \mathbf{A}_i)$, $\mathbf{A}_i \in M_{k,k}$.
Let $\mathbf{m}_{:ji}$ and $\mathbf{a}_{:ji}$ be the $j$-th column of
$\mathbf{M}_i$ and $\mathbf{A}_i$ respectively,
then $\mathbf{m}_{:ji} = \mathbf{X}_i \mathbf{a}_{:ji}$.
Because $\mathbf{X}_i$ has orthonormal columns, we have:
\begin{equation} \label{eq:measure-subspace}
  d \mathbf{m} = \prod_{i=1}^l d \mathbf{m}_i
  = \prod_{i=1}^l \prod_{j=1}^k d \mathbf{m}_{:ji}
  = \prod_{i=1}^l \prod_{j=1}^k d (\mathbf{X}_i \mathbf{a}_{:ji})
  = \prod_{i=1}^l \prod_{j=1}^k d \mathbf{a}_{:ji}
  = d \mathbf{a}
\end{equation}
Here, $\mathbf{a} = \text{vec}(\mathbf{\EuScript{A}}) \in \mathbb{R}^{kkl}$
and $\mathbf{\EuScript{A}}$ is the $k \times k \times l$ array with frontal slices $\mathbf{A}_i$.
Replacing the integration domain $S$ with $\mathbb{R}^{kk l}$, we have:
\begin{equation*}
  \circledast \propto \int_{\mathbb{R}^{kk l}} p(\mathbf{m}_* | \mathbf{m})
  ~p(\mathbf{m} | \mathfrak{X})
  d \mathbf{a}
\end{equation*}

Let $N(\mathbf{x}; \boldsymbol{\mu}, \boldsymbol{\Sigma})$ denote the value at $\mathbf{x}$
of the Gaussian PDF with mean $\boldsymbol{\mu}$ and covariance matrix $\boldsymbol{\Sigma}$.
From \cref{eq:predictive-point-subspace}, we have:
\begin{align} \label{eq:predictive-point-subspace-proof}
  p(\mathbf{m}_* | \mathbf{m})
  &= N_{nk}(\mathbf{m}_*; \mathbf{K}_{12} \mathbf{K}_{22}^{-1} \mathbf{m},
    \mathbf{I}_{nk} - \mathbf{K}_{12} \mathbf{K}_{22}^{-1} \mathbf{K}_{12}^T) \nonumber \\
  &\propto
    \exp\left(-\frac{1}{2} (\mathbf{m}_* - \mathbf{K}_{12} \mathbf{K}_{22}^{-1} \mathbf{m})^T
    \mathbf{S}^{-1} (\mathbf{m}_* - \mathbf{K}_{12} \mathbf{K}_{22}^{-1} \mathbf{m})\right)
\end{align}
Here, $\mathbf{S} = \mathbf{I}_{nk} - \mathbf{K}_{12} \mathbf{K}_{22}^{-1} \mathbf{K}_{12}^T$.
By computation rules of the Kronecker product:
\begin{align} \label{eq:S}
  \mathbf{S}
  &= \mathbf{I}_{nk} - (\mathbf{k}_l^T \otimes \mathbf{I}_{nk})
    (\mathbf{K}_l \otimes \mathbf{I}_{nk})^{-1} (\mathbf{k}_l^T \otimes \mathbf{I}_{nk})^T\nonumber\\
  &= \mathbf{I}_{nk} - (\mathbf{k}_l^T \mathbf{K}_l^{-1} \mathbf{k}_l) \otimes \mathbf{I}_{nk}
    = (1 - \mathbf{k}_l^T \mathbf{K}_l^{-1} \mathbf{k}_l) \mathbf{I}_{nk}
    = \varepsilon^2 \mathbf{I}_{nk}
\end{align}
We denote noise variance $\varepsilon^2 = 1 - \mathbf{k}_l^T \mathbf{K}_l^{-1} \mathbf{k}_l$.
Since $\mathbf{m} = (\text{vec}(\mathbf{X}_i \mathbf{A}_i))_{i=1}^l$,
by computation rules of the Kronecker product,
we can write $\mathbf{K}_{12} \mathbf{K}_{22}^{-1} \mathbf{m}$ as:
\begin{equation*}
  (\mathbf{k}_l^T \otimes \mathbf{I}_{nk}) (\mathbf{K}_l \otimes \mathbf{I}_{nk})^{-1} \mathbf{m}
  = ((\mathbf{k}_l^T \mathbf{K}_l^{-1}) \otimes \mathbf{I}_{nk}) \mathbf{m}
  = \sum_{i=1}^l (\mathbf{k}_l^T \mathbf{K}_l^{-1} \mathbf{e}_i) \text{vec}(\mathbf{X}_i \mathbf{A}_i)
\end{equation*}
Since $\text{vec}()$ is a linear operator, we have
\begin{equation*}
  \mathbf{K}_{12} \mathbf{K}_{22}^{-1} \mathbf{m} = \text{vec}\left(\sum_{i=1}^l
  (\mathbf{k}_l^T \mathbf{K}_l^{-1} \mathbf{e}_i)  \mathbf{X}_i \mathbf{A}_i\right)
  = \text{vec}\left(\sum_{i=1}^l \mathbf{X}_i (\mathbf{k}_l^T \mathbf{K}_l^{-1} \mathbf{e}_i)
  \mathbf{I}_k \mathbf{A}_i\right)
\end{equation*}
Let $\mathbf{A}_{(13\times2)}$ be the matricization of $\mathbf{\EuScript{A}}$
by combining the matrices $\mathbf{A}_i$ by rows.
Recall that $\mathbf{X}$ combines $\mathbf{X}_i$ by columns, we have
\begin{equation*}
  \mathbf{K}_{12} \mathbf{K}_{22}^{-1} \mathbf{m} = \text{vec}\left(\mathbf{X}
  (\diag(\mathbf{K}_l^{-1} \mathbf{k}_l) \otimes \mathbf{I}_k) \mathbf{A}_{(13\times2)}\right)
  = \text{vec}\left(\widetilde{\mathbf{X}} \mathbf{A}_{(13\times2)}\right)
\end{equation*}
Here,
$\widetilde{\mathbf{X}} = \mathbf{X} (\diag(\mathbf{K}_l^{-1} \mathbf{k}_l) \otimes \mathbf{I}_k)$.
By the ``vec trick'' of the Kronecker product, we have
\begin{equation} \label{eq:KKv}
  \mathbf{K}_{12} \mathbf{K}_{22}^{-1} \mathbf{m}
  = (\mathbf{I}_k \otimes \widetilde{\mathbf{X}}) \text{vec}(\mathbf{A}_{(13\times2)})
  = (\mathbf{I}_k \otimes \widetilde{\mathbf{X}}) \mathbf{a}_{(13\times2)}
\end{equation}
Here $\mathbf{a}_{(13\times2)} = \text{vec}(\mathbf{A}_{(13\times2)})$.
Substituting \cref{eq:S,eq:KKv} into \cref{eq:predictive-point-subspace-proof}, we have
\begin{equation} \label{eq:predictive-point-subspace-simplify}
  p(\mathbf{m}_* | \mathbf{m}) \propto
  \exp\left(-\frac{1}{2} \varepsilon^{-2} \|\mathbf{m}_* -
    (\mathbf{I}_k \otimes \widetilde{\mathbf{X}}) \mathbf{a}_{(13\times2)}\|^2\right)
\end{equation}

From \cref{eq:posterior-subspace}, $p(\mathbf{m} | \mathfrak{X}) \propto
\exp\{-\frac{1}{2} \mathbf{m}^T (\mathbf{K}_l^{-1} \otimes \mathbf{I}_{nk}) \mathbf{m}\}$,
where $\mathbf{m} = (\text{vec}(\mathbf{X}_i \mathbf{A}_i))_{i=1}^l$.
Note that matrix inverse and the Kronecker product commute.
Let $\bar{k}_{ij} = [\mathbf{K}_l^{-1}]_{ij}$.
Expand the Kronecker product and use properties of the trace, we have
\begin{align*}
  &\mathbf{m}^T (\mathbf{K}_l^{-1} \otimes \mathbf{I}_{nk}) \mathbf{m} \\
  =& \sum_{i=1}^l \sum_{j=1}^l \bar{k}_{ij} \text{vec}(\mathbf{X}_i \mathbf{A}_i)^T
     \text{vec}(\mathbf{X}_j \mathbf{A}_j)
     = \sum_{i=1}^l \sum_{j=1}^l \bar{k}_{ij} \text{tr}\left((\mathbf{X}_i \mathbf{A}_i)^T
     (\mathbf{X}_j \mathbf{A}_j)\right) \\
  =&~\text{tr}\left(\sum_{i=1}^l \sum_{j=1}^l \bar{k}_{ij}
     (\mathbf{X}_i \mathbf{A}_i)^T (\mathbf{X}_j \mathbf{A}_j)\right)
     = \text{tr}\left(\sum_{i=1}^l \sum_{j=1}^l (\mathbf{X}_i \mathbf{A}_i)^T
     (\bar{k}_{ij} \mathbf{I}_{n}) (\mathbf{X}_j \mathbf{A}_j)\right)
\end{align*}
Let $(\mathbf{X}_i \mathbf{A}_i)_{i=1}^l$ be the matrix
combining $\mathbf{X}_i \mathbf{A}_i$ by rows.
Let $\mathbb{X} = \diag(\mathbf{X}_i)_{i=1}^l$, then
$\mathbb{X} \mathbf{A}_{(13\times 2)} = (\mathbf{X}_i \mathbf{A}_i)_{i=1}^l$.
Reconstruct a Kronecker product, we have
\begin{align*}
  \mathbf{m}^T (\mathbf{K}_l^{-1} \otimes \mathbf{I}_{nk}) \mathbf{m}
  &= \text{tr}\left([(\mathbf{X}_i \mathbf{A}_i)_{i=1}^l]^T (\mathbf{K}_l^{-1} \otimes \mathbf{I}_{n})
    [(\mathbf{X}_j \mathbf{A}_j)_{j=1}^l]\right) \\
  &= \text{tr}\left(\mathbf{A}_{(13\times 2)}^T \mathbb{X}^T
    (\mathbf{K}_l^{-1} \otimes \mathbf{I}_{n}) \mathbb{X} \mathbf{A}_{(13\times 2)}
    \right)
\end{align*}
Let $\breve{\Box} = \mathbb{X}^T (\mathbf{K}_l^{-1} \otimes \mathbf{I}_{n}) \mathbb{X}$.
With the ``vec trick'', we have
\begin{align}
  \mathbf{m}^T (\mathbf{K}_l^{-1} \otimes \mathbf{I}_{nk}) \mathbf{m}
  &= \text{tr}\left(\mathbf{A}_{(13\times 2)}^T~\breve{\Box}~\mathbf{A}_{(13\times 2)} \right)
    = \text{vec}(\mathbf{A}_{(13\times 2)})^T \text{vec}(\breve{\Box}~\mathbf{A}_{(13\times 2)})
    \nonumber \\
  &= \text{vec}(\mathbf{A}_{(13\times 2)})^T
    (\mathbf{I}_k \otimes \breve{\Box})~\text{vec}(\mathbf{A}_{(13\times 2)})
    = \mathbf{a}_{(13\times 2)}^T (\mathbf{I}_k \otimes \breve{\Box}) \mathbf{a}_{(13\times 2)}
    \label{eq:mKm}
\end{align}
So the posterior distribution has the form:
\begin{equation} \label{eq:posterior-subspace-simplify}
  p(\mathbf{m} | \mathfrak{X}) \propto \exp\{-\frac{1}{2}
  \mathbf{a}_{(13\times 2)}^T (\mathbf{I}_k \otimes \breve{\Box}) \mathbf{a}_{(13\times 2)} \}
\end{equation}

Substitute \cref{eq:predictive-point-subspace-simplify,eq:posterior-subspace-simplify}
into $\circledast$, we have:
\begin{equation*}
  \circledast \propto \int_{\mathbb{R}^{kkl}} \exp\left(-\frac{1}{2} \left[\varepsilon^{-2}
  \|\mathbf{m}_* - (\mathbf{I}_k \otimes \widetilde{\mathbf{X}}) \mathbf{a}_{(13\times2)}\|^2
  + \mathbf{a}_{(13\times 2)}^T (\mathbf{I}_k \otimes \breve{\Box}) \mathbf{a}_{(13\times 2)}
  \right] \right) d\mathbf{a}
\end{equation*}
Note that we can expand the inner product to have:
\begin{equation*}
  \|\mathbf{m}_* - (\mathbf{I}_k \otimes \widetilde{\mathbf{X}}) \mathbf{a}_{(13\times2)}\|^2
  = \|\mathbf{m}_*\|^2 - 2 \mathbf{m}_*^T
  (\mathbf{I}_k \otimes \widetilde{\mathbf{X}}) \mathbf{a}_{(13\times2)}
  + \mathbf{a}_{(13\times2)}^T
  (\mathbf{I}_k \otimes (\widetilde{\mathbf{X}}^T \widetilde{\mathbf{X}})) \mathbf{a}_{(13\times2)}
\end{equation*}
Denote $\boldsymbol{\Sigma}_c^{-1} = \mathbf{I}_k \otimes
(\varepsilon^{-2} \widetilde{\mathbf{X}}^T\widetilde{\mathbf{X}} + \breve{\Box})$
and $\mathbf{m}_c^T \boldsymbol{\Sigma}_c^{-1} = \varepsilon^{-2} \mathbf{m}_*^T
(\mathbf{I}_k \otimes \widetilde{\mathbf{X}})$.
Because $d \mathbf{a} = d \mathbf{a}_{(13\times2)}$, we have
\begin{align*} \label{eq:integral-proof-1}
  \circledast
  &\propto
    \int_{\mathbb{R}^{kkl}} \exp\left(-\frac{1}{2} \varepsilon^{-2} \|\mathbf{m}_*\|^2
    + \mathbf{m}_c^T \boldsymbol{\Sigma}_c^{-1} \mathbf{a}_{(13\times2)}
    -\frac{1}{2} \mathbf{a}_{(13\times2)}^T \boldsymbol{\Sigma}_c^{-1} \mathbf{a}_{(13\times2)}
    \right)~d\mathbf{a}_{(13\times2)}  \nonumber \\
  & = \det(2 \pi \boldsymbol{\Sigma}_c)^{1/2} \exp\left(-\frac{1}{2} \varepsilon^{-2}
    \|\mathbf{m}_*\|^2 + \frac{1}{2} \mathbf{m}_c^T \boldsymbol{\Sigma}_c^{-1} \mathbf{m}_c\right) 
\end{align*}
With the definitions of $\boldsymbol{\Sigma}_c^{-1}$
and $\mathbf{m}_c^T \boldsymbol{\Sigma}_c^{-1}$, we have
\begin{align*}
  & \varepsilon^{-2} \|\mathbf{m}_*\|^2 - \mathbf{m}_c^T \boldsymbol{\Sigma}_c^{-1} \mathbf{m}_c 
  = \varepsilon^{-2} \|\mathbf{m}_*\|^2 -(\mathbf{m}_c^T \boldsymbol{\Sigma}_c^{-1})
    (\boldsymbol{\Sigma}_c^{-1})^{-1}
    (\mathbf{m}_c^T \boldsymbol{\Sigma}_c^{-1})^T \\
  =& \varepsilon^{-2} \|\mathbf{m}_*\|^2 - \varepsilon^{-4} \mathbf{m}_*^T
    (\mathbf{I}_k \otimes \widetilde{\mathbf{X}})  
    (\mathbf{I}_k \otimes (\varepsilon^{-2} \widetilde{\mathbf{X}}^T
    \widetilde{\mathbf{X}} + \breve{\Box}))^{-1}
    (\mathbf{I}_k \otimes \widetilde{\mathbf{X}})^T \mathbf{m}_* \\
  =& \mathbf{m}_*^T \left(\varepsilon^{-2} \mathbf{I}_{nk}
    - \varepsilon^{-4}\mathbf{I}_k \otimes (\widetilde{\mathbf{X}}
    (\varepsilon^{-2} \widetilde{\mathbf{X}}^T \widetilde{\mathbf{X}} + \breve{\Box})^{-1}
    \widetilde{\mathbf{X}}^T)\right) \mathbf{m}_*
  = \mathbf{m}_*^T \left(\mathbf{I}_k \otimes \boldsymbol{\Sigma}^{-1} \right) \mathbf{m}_*
\end{align*}
Here we define
\begin{equation} \label{eq:pred-covar-raw}
  \boldsymbol{\Sigma}^{-1} = \varepsilon^{-2} \mathbf{I}_n
  - \varepsilon^{-4} \widetilde{\mathbf{X}}
  (\varepsilon^{-2} \widetilde{\mathbf{X}}^T \widetilde{\mathbf{X}} + \breve{\Box})^{-1}
  \widetilde{\mathbf{X}}^T
\end{equation}
Because $\boldsymbol{\Sigma}_c$ does not depend on $\mathbf{m}_*$
but $\circledast = p(\mathbf{m}_* | \mathfrak{X})$, we have
\begin{equation} \label{eq:integral-proof-2}
  \circledast \propto \exp\left(-\frac{1}{2} \mathbf{m}_*^T
  \left(\mathbf{I}_k \otimes \boldsymbol{\Sigma}^{-1} \right) \mathbf{m}_* \right)
\end{equation}
This means that the predictive distribution is
\begin{equation} \label{eq:predictive-subspace-proof}
  \circledast := p(\mathbf{m}_* | \mathfrak{X})
  = N_{nk}(\mathbf{m}_*; 0, \mathbf{I}_k \otimes \boldsymbol{\Sigma})
\end{equation}

Now we simplify $\boldsymbol{\Sigma}$.
Recall that $\widetilde{\mathbf{X}} =
\mathbf{X} (\diag(\mathbf{K}_l^{-1} \mathbf{k}_l) \otimes \mathbf{I}_k)$.
Let $\mathbf{v} = \mathbf{K}_l^{-1} \mathbf{k}_l$.
Using the definition and properties of the Kronecker product, we have the following:
\begin{gather*}
  \widetilde{\mathbf{X}} = \mathbf{X} (\diag(\mathbf{v}) \otimes \mathbf{I}_k) =
  (\mathbf{v} \otimes \mathbf{I}_n) \mathbb{X} \\
  \widetilde{\mathbf{X}}^T \widetilde{\mathbf{X}}= 
  \mathbb{X}^T (\mathbf{v} \otimes \mathbf{I}_n)^T (\mathbf{v} \otimes \mathbf{I}_n) \mathbb{X} =
  \mathbb{X}^T [(\mathbf{v}^T \mathbf{v}) \otimes \mathbf{I}_n)] \mathbb{X}
\end{gather*}
Recall that $\breve{\Box} = \mathbb{X}^T (\mathbf{K}_l^{-1} \otimes \mathbf{I}_{n}) \mathbb{X}$,
from \cref{eq:pred-covar-raw} and the above, we have
\begin{align*} \boldsymbol{\Sigma}^{-1}
  &= \varepsilon^{-2} \mathbf{I}_n - \varepsilon^{-4} \widetilde{\mathbf{X}} 
    \{\varepsilon^{-2} \mathbb{X}^T [(\mathbf{v}^T \mathbf{v}) \otimes \mathbf{I}_n)] \mathbb{X}
    + \mathbb{X}^T (\mathbf{K}_l^{-1} \otimes \mathbf{I}_{n}) \mathbb{X}\}^{-1}
    \widetilde{\mathbf{X}}^T \\
  &= \varepsilon^{-2} \mathbf{I}_n - \varepsilon^{-4} \widetilde{\mathbf{X}} 
    \{\mathbb{X}^T [(\varepsilon^{-2} \mathbf{v}^T \mathbf{v} + \mathbf{K}_l^{-1})
    \otimes \mathbf{I}_n)] \mathbb{X}\}^{-1} \widetilde{\mathbf{X}}^T
\end{align*}
Let $\mathbf{D}_{\mathbf{v}} = \diag(\mathbf{v})$,
then $\widetilde{\mathbf{X}} = \mathbf{X} (\mathbf{D}_{\mathbf{v}} \otimes \mathbf{I}_k)$.
Since $\mathbb{X} (\mathbf{D}_{\mathbf{v}}^{-1} \otimes \mathbf{I}_k) =
(\mathbf{D}_{\mathbf{v}}^{-1} \otimes \mathbf{I}_n) \mathbb{X}$, we have
\begin{align*} \boldsymbol{\Sigma}^{-1}
  &= \varepsilon^{-2} \mathbf{I}_n - \varepsilon^{-4} \mathbf{X}
    (\mathbf{D}_{\mathbf{v}} \otimes \mathbf{I}_k) \{\mathbb{X}^T
    [(\varepsilon^{-2} \mathbf{v}^T \mathbf{v} + \mathbf{K}_l^{-1}) \otimes \mathbf{I}_n)]
    \mathbb{X}\}^{-1} (\mathbf{D}_{\mathbf{v}} \otimes \mathbf{I}_k) \mathbf{X}^T \\
  &= \varepsilon^{-2} \mathbf{I}_n - \varepsilon^{-4} \mathbf{X}
    \{(\mathbf{D}_{\mathbf{v}} \otimes \mathbf{I}_k)^{-1} \mathbb{X}^T
    [(\varepsilon^{-2} \mathbf{v}^T \mathbf{v} + \mathbf{K}_l^{-1}) \otimes \mathbf{I}_n)]
    \mathbb{X} (\mathbf{D}_{\mathbf{v}} \otimes \mathbf{I}_k)^{-1}\}^{-1} \mathbf{X}^T \\
  &= \varepsilon^{-2} \mathbf{I}_n - \varepsilon^{-4} \mathbf{X}
    \{\mathbb{X}^T (\mathbf{D}_{\mathbf{v}}^{-1} \otimes \mathbf{I}_n)
    [(\varepsilon^{-2} \mathbf{v}^T \mathbf{v} + \mathbf{K}_l^{-1}) \otimes \mathbf{I}_n)]
    (\mathbf{D}_{\mathbf{v}}^{-1} \otimes \mathbf{I}_n) \mathbb{X}\}^{-1} \mathbf{X}^T \\
  &= \varepsilon^{-2} \mathbf{I}_n - \varepsilon^{-4} \mathbf{X} \{\mathbb{X}^T [(\varepsilon^{-2}
    \mathbf{D}_{\mathbf{v}}^{-1} \mathbf{v}^T \mathbf{v} \mathbf{D}_{\mathbf{v}}^{-1}
    + \mathbf{D}_{\mathbf{v}}^{-1} \mathbf{K}_l^{-1} \mathbf{D}_{\mathbf{v}}^{-1})
    \otimes \mathbf{I}_n)] \mathbb{X}\}^{-1} \mathbf{X}^T \\
  &= \varepsilon^{-2} \mathbf{I}_n - \varepsilon^{-4} \mathbf{X} \{\mathbb{X}^T
    [(\varepsilon^{-2} \mathbf{1}_l \mathbf{1}_l^T
    + \mathbf{D}_{\mathbf{v}}^{-1} \mathbf{K}_l^{-1} \mathbf{D}_{\mathbf{v}}^{-1})
    \otimes \mathbf{I}_n)] \mathbb{X}\}^{-1} \mathbf{X}^T \\
  &= \varepsilon^{-2} \mathbf{I}_n - \varepsilon^{-4} \mathbf{X} [\mathbb{X}^T
    (\boldsymbol{\Omega} \otimes \mathbf{I}_n) \mathbb{X}]^{-1} \mathbf{X}^T
\end{align*}
In the last step, we define $\boldsymbol{\Omega}= \varepsilon^{-2} \mathbf{1}_l \mathbf{1}_l^T
+ \mathbf{D}_{\mathbf{v}}^{-1} \mathbf{K}_l^{-1} \mathbf{D}_{\mathbf{v}}^{-1}$.
With the Woodbury identity:
\begin{equation*}
  (\mathbf{A} + \mathbf{C} \mathbf{B} \mathbf{C}^T)^{-1} =
  \mathbf{A}^{-1} - \mathbf{A}^{-1} \mathbf{C}(\mathbf{B}^{-1} +
  \mathbf{C}^T \mathbf{A}^{-1} \mathbf{C})^{-1} \mathbf{C}^T \mathbf{A}^{-1}
\end{equation*}
we substitute $\mathbf{A} = \varepsilon^{-2} \mathbf{I}_n$,
$\mathbf{B} = -[\mathbb{X}^T (\boldsymbol{\Omega} \otimes \mathbf{I}_n) \mathbb{X}]^{-1}$,
and $\mathbf{C} = \varepsilon^{-2} \mathbf{X}$. This gives:
\begin{equation*}
  \boldsymbol{\Sigma} = \varepsilon^2 \mathbf{I}_n - \mathbf{X}
  [-\mathbb{X}^T (\boldsymbol{\Omega} \otimes \mathbf{I}_n) \mathbb{X}
  + \varepsilon^{-2} \mathbf{X}^T \mathbf{X}]^{-1} \mathbf{X}^T
\end{equation*}
Note that $\mathbf{X} = (\mathbf{1}_l^T \otimes \mathbf{I}_n) \mathbb{X}$, so we have:
\begin{equation*}
  \mathbf{X}^T \mathbf{X} = \mathbb{X}^T (\mathbf{1}_l^T \otimes \mathbf{I}_n)^T
  (\mathbf{1}_l^T \otimes \mathbf{I}_n) \mathbb{X}
  = \mathbb{X}^T [(\mathbf{1}_l \mathbf{1}_l^T) \otimes \mathbf{I}_n] \mathbb{X}
\end{equation*}
Let $\widetilde{\mathbf{K}}_l = (\mathbf{D}_{\mathbf{v}} \mathbf{K}_l \mathbf{D}_{\mathbf{v}})^{-1}$,
then $\boldsymbol{\Omega}= \varepsilon^{-2} \mathbf{1}_l \mathbf{1}_l^T + \widetilde{\mathbf{K}}_l$.
We have:
\begin{align}
  \boldsymbol{\Sigma}
  &= \varepsilon^2 \mathbf{I}_n - \mathbf{X}
    [-\mathbb{X}^T (\boldsymbol{\Omega} \otimes \mathbf{I}_n) \mathbb{X}
    + \varepsilon^{-2} \mathbb{X}^T [(\mathbf{1}_l \mathbf{1}_l^T) \otimes \mathbf{I}_n]
    \mathbb{X}]^{-1} \mathbf{X}^T \nonumber \\
  &= \varepsilon^2 \mathbf{I}_n + \mathbf{X}
    \{\mathbb{X}^T [(\boldsymbol{\Omega} - \varepsilon^{-2} \mathbf{1}_l \mathbf{1}_l^T)
    \otimes \mathbf{I}_n] \mathbb{X}\}^{-1} \mathbf{X}^T \nonumber \\
  &= \varepsilon^2 \mathbf{I}_n + \mathbf{X}
    [\mathbb{X}^T (\widetilde{\mathbf{K}}_l \otimes \mathbf{I}_n) \mathbb{X}]^{-1} \mathbf{X}^T
    \label{eq:covariance-subspace-proof}
\end{align}

With \cref{eq:predictive-subspace-proof,eq:covariance-subspace-proof}, we complete the proof.

\section{Joint distributions and random functions on Grassmann manifold}
\label{sup:joint-distribution}

In the main text we focus on point predictions on the Grassmann manifold,
which is enough for PROM purposes.
But more generally, our GP model induces a family of joint distributions on Grassmann manifolds,
and can be used to generate random subspace-valued functions.
Neither of these problems have been explored in the literature.

From \cref{sec:subspace-regression}, we see that
for any finite collection of parameter points $\boldsymbol{\theta} = (\boldsymbol{\theta}_i)_{i=1}^l$,
our GP model gives a collection of random points on the Grassmann manifold
$\mathfrak{M}_i = \text{span}(\text{vec}^{-1}(\bar{f}(\boldsymbol{\theta}_i)))$,
whose marginal distributions are uniform: $\mathfrak{M}_i \sim \text{Uniform}(G_{k,n})$.
For each $i \in \{2, \cdots, l\}$,
let $\boldsymbol{\Sigma}_{\le i}$ be defined by
$\boldsymbol{\theta}_{\le i} = (\boldsymbol{\theta}_j)_{j=1}^i$
and $\mathfrak{M}_{< i} = (\mathfrak{M}_j)_{j=1}^{i-1}$
as in \cref{eq:covariance-subspace}.
Then we have conditional distributions
$\mathfrak{M}_i | \mathfrak{M}_{< i} \sim \text{MACG}(\boldsymbol{\Sigma}_{\le i})$.
Combining the marginal and conditional distributions,
we have a joint distribution on the Grassmann manifold,
parameterized by $\boldsymbol{\theta}$:
\begin{equation} \label{eq:joint-distribution}
  (\mathfrak{M}_i)_{i=1}^l \sim \text{Uniform}(G_{k,n})
  \prod_{i=2}^l \text{MACG}(\boldsymbol{\Sigma}_{\le i})
\end{equation}

GPS can be used to generate random subspace-valued functions.
Suppose that $\boldsymbol{\theta}$ is a sample grid to evaluate the random function,
then we can use \cref{eq:joint-distribution} to generate a sample path sequentially.
The method to sample $\text{MACG}(\boldsymbol{\Sigma})$, including the uniform distribution,
is implied in \cref{sub:Grassmann}, which requires $\boldsymbol{\Sigma}^{1/2}$.
If we compute the EVD of $\boldsymbol{\Sigma}$ as in \cref{sec:prediction},
then we have $\boldsymbol{\Sigma}^{1/2} =
\mathbf{V} \diag(\sqrt{\sigma_i^2 + \varepsilon^2})_{i=1}^r \mathbf{V}^T + \varepsilon \mathbf{I}_n$.
We summarize the overall sampling procedure in \cref{alg:GPS-sample}.

\alglanguage{pseudocode}
\begin{algorithm}[h]
  \caption{GPS: Sampling}
  \label{alg:GPS-sample}
  \begin{algorithmic}[1] %
    \Require correlation function $k(\cdot, \cdot)$.
    \Input sample grid $(\boldsymbol{\theta}_i)_{i=1}^l$.
    \State Generate random matrix: $\mathbf{Z} \in M_{n,k}$, $z_{ij} \sim N(0, 1)$.
    \State Orthonormalization: $\mathbf{X}_1 \gets \pi(\mathbf{Z})$.
    \For{$i$ in $2, \cdots, l$}
    \State Generate random matrix: $\mathbf{Z} \in M_{n,k}$, $z_{ij} \sim N(0, 1)$.
    \State %
    Run \cref{alg:GPS-preprocess,alg:GPS-predict-EVD} with arguments
    $\mathbf{X}_{< i}$ and $(\boldsymbol{\theta}_{< i}, \boldsymbol{\theta}_i, r)$.
    \State Matrix multiplication: $\mathbf{M} \gets \mathbf{V}
    \diag(\sqrt{\mathring{\boldsymbol{\lambda}} + \varepsilon^2} - \varepsilon) \mathbf{V}^T
    \mathbf{Z} + \varepsilon \mathbf{Z}$
    \State Orthonormalization: $\mathbf{X}_i \gets \pi(\mathbf{M})$.
    \EndFor
    \Output Stiefel representations of subspaces $(\mathbf{X}_i)_{i=1}^l$.
    \Note Projection $\pi(\mathbf{M}) = \mathbf{U} \mathbf{W}^T$, where
    $\mathbf{M} = \mathbf{U} \diag(\sigma) \mathbf{W}^T$ is a thin SVD.
  \end{algorithmic}
\end{algorithm}

\section{Gradient of LOOCV predictive error}
\label{sup:gradient}

The gradient of the LOOCV predictive error can also be computed.
Denote $d_i = d_g(\mathbf{X}_i, \mathbf{V}_{-i})$
and let $\partial$ denote the partial derivative with respect to a scalar hyperparameter.
With \cref{eq:LOOCV-error} and chain rule, we have:
\begin{equation} \label{eq:d-e2}
  \partial \epsilon_2 = \sum_{i=1}^l \partial d_i^2
  = -2 \sum_{i=1}^l \sum_{j=1}^k (\arccos \sigma_j) (1-\sigma_j^2)^{-1/2} \partial \sigma_j
\end{equation}
Here, $\sigma_j = \sigma_j(\mathbf{X}_i^T \mathbf{V}_{-i})
= \sigma_j(\widetilde{\mathbf{C}}_i^T \mathring{\mathbf{V}}_{-i})$.
Let $\widetilde{\mathbf{C}}_i^T \mathring{\mathbf{V}}_{-i}
= \widehat{\mathbf{V}} \diag(\boldsymbol{\sigma}) \widehat{\mathbf{W}}^T$ be a thin SVD.
Using the derivative of a singular value, see for example \cite[p. 170]{Strang2019},
we have:
\begin{equation} \label{eq:d-sigma}
  \partial \sigma_j
  = \widehat{\mathbf{v}}_j^T (\widetilde{\mathbf{C}}_i^T \partial \mathring{\mathbf{V}})
  \widehat{\mathbf{w}}_j
  = \widehat{\mathbf{v}}_j^T \widetilde{\mathbf{C}}_i^T (\partial \mathring{\mathbf{V}})
  \widehat{\mathbf{w}}_j
\end{equation}
Recall that $\mathring{\mathbf{V}}$ consists of the top-$k$ eigenvectors of $\mathbf{S}_{-i}$.
Let $(\mathring{\lambda}_p, \mathring{\mathbf{v}}_p)$
be the $p$-th eigenpair of $\mathbf{S}_{-i}$, $p = 1, \cdots, k$.
Using the derivative of an eigenvector of a symmetric matrix,
see for example \cite[Thm 8.9]{Magnus2019}, we have:
\begin{equation} \label{eq:d-v}
  \partial \mathring{\mathbf{v}}_p = (\mathring{\lambda}_p \mathbf{I} - \mathbf{S}_{-i})^\dagger
  (\partial \mathbf{S}_{-i}) \mathring{\mathbf{v}}_p
\end{equation}
Here, $\dagger$ denotes the Moore--Penrose inverse.
Let $\mathbf{S}_{-i} = \mathring{\mathbf{Q}} \diag(\mathring{\boldsymbol{\lambda}})
\mathring{\mathbf{Q}}^T$ be an EVD, then we have
$(\mathring{\lambda}_p \mathbf{I} - \mathbf{S}_{-i})^\dagger = \mathring{\mathbf{V}}
\diag\{(\mathring{\lambda}_p - \mathring{\lambda}_q)^{-1}\}_{q=1}^r \mathring{\mathbf{V}}^T$.
Recall that $\mathbf{S}_{-i} =
\widetilde{\mathbf{C}}_{-i} (\boldsymbol{\Pi}_{-i})^{-1} \widetilde{\mathbf{C}}_{-i}^T$, we have:
\begin{equation} \label{eq:d-S}
  \partial \mathbf{S}_{-i} = - \widetilde{\mathbf{C}}_{-i} (\boldsymbol{\Pi}_{-i})^{-1}
  (\partial \boldsymbol{\Pi}_{-i}) (\boldsymbol{\Pi}_{-i})^{-1} \widetilde{\mathbf{C}}_{-i}^T
\end{equation}
Recall that $\boldsymbol{\Pi}_{-i} = \Box_{-i} \circ (\Delta_{-i} \otimes \mathbf{J}_k)$,
$\Delta_{-i} = [\bar{k}_{pq} \bar{k}_{ii} / (\bar{k}_{ip} \bar{k}_{iq}) - 1]_{p,q \ne i}$,
and $\bar{\mathbf{K}} = \mathbf{K}^{-1}$, we have:
\begin{gather}
  \partial \boldsymbol{\Pi}_{-i} = \Box_{-i} \circ [(\partial \Delta_{-i}) \otimes \mathbf{J}_k]
  \nonumber \\
  \partial [\Delta_{-i}]_{pq} = ([\Delta_{-i}]_{pq} + 1)
  (\bar{k}_{pq}^{-1} \partial \bar{k}_{pq} + \bar{k}_{ii}^{-1} \partial \bar{k}_{ii}
  - \bar{k}_{ip}^{-1} \partial \bar{k}_{ip} - \bar{k}_{iq}^{-1} \partial \bar{k}_{iq}) \label{eq:d-Pi}\\
  \partial \bar{k}_{pq} = [\partial \mathbf{K}^{-1}]_{pq} =
  [-\mathbf{K}^{-1} (\partial \mathbf{K}) \mathbf{K}^{-1}]_{pq} \nonumber
\end{gather}
Combining \cref{eq:d-e2,eq:d-sigma,eq:d-v,eq:d-S,eq:d-Pi},
we can compute the partial derivative $\partial \epsilon_2$
of the LOOCV predictive error with respect to any hyperparameter,
as long as we can compute the partial derivative $\partial k$ of the correlation function.
For the SE kernel in \cref{eq:SE-kernel} for example,
$\partial k / \partial \beta_i = (\theta_i-\theta_i')^2 \beta_i^{-3} k$.
We omit a formal algorithm for the gradient computation,
since it is straightforward given these equations.

We point out one way to speed up the evaluation of \cref{eq:d-v,eq:d-S}.
Because the eigenvalues of $\mathbf{S}_{-i}$ decline rapidly, we have
\begin{equation}
  (\mathring{\lambda}_p \mathbf{I} - \mathbf{S}_{-i})^\dagger \approx
  \mathring{\mathbf{V}} \diag\{(\mathring{\lambda}_p - \mathring{\lambda}_q)^{-1}\}_{q=1}^\tau
  \mathring{\mathbf{V}}^T
  + \mathring{\lambda}_p^{-1} (\mathbf{I} - \mathring{\mathbf{V}} \mathring{\mathbf{V}}^T)
\end{equation}
This approximation is accurate for any $p \in \{1, \cdots, k\}$,
as long as $\tau - k$ is reasonably large; for example, we can set $\tau = 2 k$.
To compute the approximation we only need the top $\tau$ eigenpairs of $\mathbf{S}_{-i}$.
Since $r \approx kl > 2k$, the truncated EVD can be substantially faster than a full EVD.
\Cref{alg:approx-d-v} gives an efficient procedure to compute $\partial\mathring{\mathbf{v}}_p$
approximately given $\mathring{\mathbf{v}}_p$ and $\partial \boldsymbol{\Pi}_{-i}$.

\alglanguage{pseudocode}
\begin{algorithm}[h]
  \caption{Approximate Computation of Derivative of an Eigenvector}
  \label{alg:approx-d-v}
  \begin{algorithmic}[1] %
    \Note This procedure evaluates $\partial \mathring{\mathbf{v}}_p$ via \cref{eq:d-v,eq:d-S}
    given $(\mathring{\mathbf{v}}_p, \partial \boldsymbol{\Pi}_{-i})$.
    \Require
    $(\mathbf{L}, \widetilde{\mathbf{L}}, \mathring{\mathbf{V}}, \mathring{\boldsymbol{\lambda}})$
    from \Cref{alg:GPS-LOOCV-error}.
    \State $\mathbf{v} \gets \text{solve}(\mathbf{L}^T,
    \widetilde{\mathbf{L}} \mathring{\mathbf{v}}_p)$
    \State $\mathbf{v} \gets \text{solve}(\mathbf{L}, (\partial \boldsymbol{\Pi}_{-i}) \mathbf{v})$
    \State $\mathbf{u} \gets -\mathring{\mathbf{V}}^T (\widetilde{\mathbf{L}} \mathbf{v})$
    \State $\mathbf{w} \gets \diag\left\{(\mathring{\lambda}_p - \mathring{\lambda}_q)^{-1}
      - \mathring{\lambda}_p^{-1}\right\}_{q=1}^\tau$
    \State $\partial\mathring{\mathbf{v}}_p \gets \mathring{\mathbf{V}} \mathbf{w}
    + \mathring{\lambda}_p^{-1} \mathbf{v}$
  \end{algorithmic}
\end{algorithm}

If the gradient is computed along with the LOOCV error,
the additional cost is dominated by (1) the extended truncated EVD of $\mathbf{S}_{-i}$ for $l$ times
and (2) the evaluation of \cref{alg:approx-d-v} for $kl$ times.
Since the additional cost of truncated EVD takes about $\mathcal{O}(k^2 l^2 (\tau - k))$ flops,
with $\tau = 2k$, part (1) takes about $\mathcal{O}(k^3 l^3)$ flops.
Since \cref{alg:approx-d-v} takes about $12 k^2 l^2$ flops, part (2) takes about $12 k^3 l^3$ flops.
The overall additional cost is about $12 k^3 l^3 + \mathcal{O}(k^3 l^3)$ flops per gradient evaluation,
where the coefficient of the second term is determined by the truncated EVD algorithm.
Compared with the $k^3 l^4$ flops for LOOCV error evaluation,
the additional cost is at a similar level, depending on $l$.

\section{Other model selection criteria}
\label{sup:criteria}

There are other model selection criteria for GP models in general.
One popular possibility is to choose the hyperparameters to maximize the marginal likelihood
with the GP integrated out.
However, this approach is less robust to model and prior misspecification than CV.
Another useful criteria is the LOOCV predictive probability density.
We derived the analytical forms of both criteria for our model,
and tried them for the numerical examples in this paper.
In all cases, the marginal likelihood prefers infinite length-scales,
inducing a singular covariance matrix.
While the LOOCV predictive probability density
can select a good length-scale for the visualization problem in \cref{sub:viz},
it also prefers infinite length-scales in other problems, probably because $n \gg k$.
We explain such behavior in this section.

The marginal likelihood of data is defined as the likelihood of data integrated over the prior.
Recall that $\mathbf{x} = (\mathbf{x}_i)_{i=1}^l$, $\mathbf{x}_i = \text{vec}(\mathbf{X}_i)$,
$\mathbf{X}_i \in V_{k,n}$, $\mathbf{m} = (\mathbf{m}_i)_{i=1}^l$,
and $\mathbf{m}_i \in \mathbb{R}^{nk}$.
Let $\mathfrak{M} = (\mathfrak{M}_i)_{i=1}^l$ and $\mathfrak{M}_i = \text{span}(\mathbf{m}_i)$,
we can write the marginal likelihood as:
\begin{equation} \label{eq:marginal-likelihood}
  p(\mathbf{x}) = \int_{\mathbb{R}^{nkl}} p(\mathbf{m}) L(\mathbf{x} | \mathfrak{M})~d \mathbf{m}
\end{equation}
But from \cref{eq:likelihood-subspace}
we have likelihood $L(\mathbf{x}_i | \mathfrak{M}_i) = 1(\mathbf{x}_i \in [\mathbf{m}_i])
= 1(\mathbf{m}_i \in [\mathbf{x}_i])$,
so the integrant in \cref{eq:marginal-likelihood} only takes positive values
for $\mathbf{m} \in \prod_{i=1}^l [\mathbf{x}_i]$,
which is a measure-zero subset of the integration domain $\mathbb{R}^{nkl}$.
This means that the marginal likelihood is identically zero.

Alternatively, we may modify the definition of marginal likelihood
to only integrate over the support $S$ of a singular likelihood,
and define a modified marginal likelihood as:
\begin{equation} \label{eq:modified-marginal-likelihood}
  \tilde{p}(\mathbf{x}) = \int_S p(\mathbf{m}) L(\mathbf{x} | \mathfrak{M})~d \mathbf{m}
\end{equation}

\begin{proposition} \label{prop:marginal-likelihood-subspace}
  Let $\breve{\Box} = \mathbb{X}^T (\mathbf{K}_l^{-1} \otimes \mathbf{I}_{n}) \mathbb{X}$.
  The log modified marginal likelihood of data is:
  \begin{equation} \label{eq:marginal-likelihood-subspace}
    \log \tilde{p}(\mathbf{x}) = -\frac{1}{2}(n-k)kl \log(2 \pi)
    - \frac{k}{2} (n \log |\mathbf{K}_l| + \log |\breve{\Box}|)
  \end{equation}
\end{proposition}

\begin{proof}[Proof of \cref{prop:marginal-likelihood-subspace}]
  As in the proof of \cref{thm:predictive-subspace},
  the support of the likelihood can be written as $S = \prod_{i=1}^l \mathfrak{X}_i^k$,
  a linear subspace of $\mathbb{R}^{nkl}$
  where $\prod_{i=1}^l [\mathbf{x}_i]$ is a full-measure subset.
  Substituting prior joint distribution
  $\mathbf{m} \sim N_{nkl}(0, \mathbf{K}_l \otimes \mathbf{I}_{nk})$
  into \cref{eq:modified-marginal-likelihood}, we have:
  \begin{align*}
    \tilde{p}(\mathbf{x}) =
    \int_{S} N_{nkl}(\mathbf{m}; 0, \mathbf{K}_l \otimes \mathbf{I}_{nk})
    \prod_{i=1}^l 1(\mathbf{m}_i \in [\mathbf{x}_i]) ~d \mathbf{m}
  \end{align*}
  With the same reasoning that leads to \cref{eq:measure-subspace},
  let $\mathbf{m}_i = \text{vec}(\mathbf{X}_i \mathbf{A}_i)$,
  then we can change the integration domain to $\mathbb{R}^{kkl}$ and
  replace $d \mathbf{m}$ with $d \mathbf{a}$, which gives:
  \begin{align*}
    \tilde{p}(\mathbf{x})
    &= \int_{\mathbb{R}^{kkl}} N_{nkl}(\mathbf{m}; 0, \mathbf{K}_l \otimes \mathbf{I}_{nk})
      ~d \mathbf{a} \\
    &= \int_{\mathbb{R}^{kkl}} \det(2 \pi \mathbf{K}_l \otimes \mathbf{I}_{kn})^{-1/2}
      \exp\left(-\frac{1}{2} \mathbf{m}^T (\mathbf{K}_l^{-1} \otimes \mathbf{I}_{kn}) \mathbf{m}\right)
      ~d \mathbf{a}
  \end{align*}
  With \cref{eq:mKm},
  let $\breve{\Box} = \mathbb{X}^T (\mathbf{K}_l^{-1} \otimes \mathbf{I}_{n}) \mathbb{X}$
  and because $d \mathbf{a} = d \mathbf{a}_{(13\times 2)}$, we have:
  \begin{align*}
    \tilde{p}(\mathbf{x})
    &= \int_{\mathbb{R}^{kkl}} \det(2 \pi \mathbf{K}_l \otimes \mathbf{I}_{kn})^{-1/2}
      \exp\left(-\frac{1}{2} \mathbf{a}_{(13\times 2)}^T (\mathbf{I}_k \otimes \breve{\Box})
      \mathbf{a}_{(13\times 2)} \right) ~d \mathbf{a}_{(13\times 2)}
  \end{align*}
  With Gaussian integral
  $\int_{\mathbb{R}^n} \exp\left(-\frac{1}{2} (\mathbf{x} - \boldsymbol{\mu})^T
    \boldsymbol{\Sigma}^{-1} (\mathbf{x} - \boldsymbol{\mu}) \right)~d\mathbf{x}
  = \det(2 \pi \boldsymbol{\Sigma})^{1/2}$, we have:
  \begin{align*}
    \tilde{p}(\mathbf{x})
    &= \det(2 \pi \mathbf{K}_l \otimes \mathbf{I}_{kn})^{-1/2}
      \det(2 \pi (\mathbf{I}_k \otimes \breve{\Box})^{-1})^{1/2} \\
    &= (2\pi)^{-nkl/2} \det(\mathbf{K}_l)^{-nk/2} (2 \pi)^{kkl/2} \det(\breve{\Box})^{-k/2} \\
    &= (2\pi)^{-(n-k)kl/2} \det(\mathbf{K}_l)^{-nk/2} \det(\breve{\Box})^{-k/2}
  \end{align*}
  Taking a logarithm gives the result in \cref{eq:marginal-likelihood-subspace}.
\end{proof}

\begin{proposition} \label{prop:marginal-likelihood-subspace-prefers-singularity}
  Maximizing the modified marginal likelihood $\tilde{p}(\mathbf{x})$
  leads to a singular covariance matrix $\mathbf{K}_l$.
\end{proposition}

\begin{proof}[Proof of \cref{prop:marginal-likelihood-subspace-prefers-singularity}]
  With \cref{prop:marginal-likelihood-subspace}, we have
  \begin{align*}
    -\log \tilde{p}(\mathbf{x}) \propto h(\boldsymbol{\beta}) :=
    n \log |\mathbf{K}_l| + \log |\mathbb{X}^T (\mathbf{K}_l^{-1} \otimes \mathbf{I}_{n}) \mathbb{X}|
  \end{align*}
  Maximizing $\tilde{p}(\mathbf{x})$ is equivalent to
  minimizing the objective function $h(\boldsymbol{\beta})$.
  Let $\mathbf{Q} = (\mathbb{X}, \mathbb{X}_{\perp})$ be an orthogonal completion of $\mathbb{X}$, then
  $|\mathbf{K}_l|^n = |\mathbf{K}_l \otimes \mathbf{I}_n| = |\mathbf{K}_l^{-1} \otimes \mathbf{I}_n|^{-1}
  = |\mathbf{Q} (\mathbf{K}_l^{-1} \otimes \mathbf{I}_n) \mathbf{Q}|^{-1}$.
  Let $\mathbf{B} = \mathbf{Q}^T (\mathbf{K}_l^{-1} \otimes \mathbf{I}_n) \mathbf{Q}$, with
  block structure $\mathbf{B} = [\mathbf{B}_{11}~\mathbf{B}_{12}; \mathbf{B}_{12}^T~\mathbf{B}_{22}]$
  where $\mathbf{B}_{11}$ is order-$kl$, then we have:
  \begin{align*}
    h(\boldsymbol{\beta})
    &= \log \frac{|\mathbb{X}^T (\mathbf{K}_l^{-1} \otimes \mathbf{I}_n) \mathbb{X}|}
      {|\mathbf{Q}^T (\mathbf{K}_l^{-1} \otimes \mathbf{I}_n) \mathbf{Q}|}
      = \log \frac{|\mathbf{B}_{11}|}{|\mathbf{B}|}
  \end{align*}
  Note that $\mathbf{B}$ is positive semi-definite and so is $\mathbf{B}_{11}$.
  By the determinant properties of a block matrix, we have
  $|\mathbf{B}| = |\mathbf{B}_{11}| |\mathbf{C}_2|$, where
  $\mathbf{C}_2 = \mathbf{B}_{22} - \mathbf{B}_{12}^T \mathbf{B}_{11}^{-1} \mathbf{B}_{12}$.
  By the inverse properties of a block matrix,
  $\mathbf{C}_2^{-1}$ is the trailing principal submatrix of
  $\mathbf{B}^{-1} = \mathbf{Q}^T (\mathbf{K}_l \otimes \mathbf{I}_n) \mathbf{Q}$.
  Therefore,
  \begin{align*}
    h(\boldsymbol{\beta}) = \log (|\mathbf{C}_2|^{-1}) = \log |\mathbf{C}_2^{-1}| =
    \log |\mathbb{X}_{\perp}^T (\mathbf{K}_l \otimes \mathbf{I}_n) \mathbb{X}_{\perp}|
  \end{align*}
  As $\mathbf{K}_l$ tends to singularity, so does
  $|\mathbb{X}_{\perp}^T (\mathbf{K}_l \otimes \mathbf{I}_n) \mathbb{X}_{\perp}|$,
  which means the objective function $h(\boldsymbol{\beta})$ drops to negative infinity.
  Therefore, minimizing $h(\boldsymbol{\beta})$ selects a singular $\mathbf{K}_l$.
\end{proof}

With an SE kernel, increasing length-scales drives $\mathbf{K}_l$ to singularity.
By \cref{prop:marginal-likelihood-subspace-prefers-singularity},
maximizing the modified marginal likelihood gives infinite length-scales.

Another model selection criteria is the log LOOCV predictive probability density.
Because the predictive distribution of our GP model is $\text{MACG}(\boldsymbol{\Sigma})$,
we have:
\begin{align*}
  \log p_\text{LOO} = \sum_{i=1}^l \log p_\text{MACG}(\mathbf{X}_i; \boldsymbol{\Sigma}_{-i})
  = -\frac{1}{2} \sum_{i=1}^l \left(k \log |\boldsymbol{\Sigma}_{-i}|
  + n \log |\mathbf{X}_i^T (\boldsymbol{\Sigma}_{-i})^{-1} \mathbf{X}_i| \right)
\end{align*}
Here, $\boldsymbol{\Sigma}_{-i}$ is defined similarly as in \cref{eq:covariance-subspace},
predicting the $i$-th sample point using the other points.
Similar to the proof of \cref{prop:marginal-likelihood-subspace-prefers-singularity},
let $\mathbf{Q}_i = (\mathbf{X}_i, \mathbf{X}_{i\perp})$ be an orthogonal completion of $\mathbf{X}_i$,
let $\mathbf{B} = \mathbf{Q}_i^T (\boldsymbol{\Sigma}_{-i})^{-1} \mathbf{Q}_i$,
and let $\mathbf{B}_{11}$ be its leading principal submatrix of order $k$, then
\begin{align*}
  \log p_\text{LOO} = -\frac{1}{2} \sum_{i=1}^l \log \frac{|\mathbf{B}_{11}|^n}{|\mathbf{B}|^k}
\end{align*}
Note that both $\mathbf{B}$ and $\mathbf{B}_{11}$ are positive semi-definite,
and of orders $n$ and $k$ respectively.
As length-scale increases, both determinants increase.
When $n$ is not way larger than $k$, as in our visualization example on $G_{1,2}$,
the LOOCV predictive probability density can select a good length-scale.
But when $n$ is much larger than $k$, as in our example PROM problems,
the numerator is less influential than the denominator,
and maximizing $p_\text{LOO}$ gives infinite length-scales.
\section{On approximating local IRKA bases}
\label{sup:irka}

\begin{figure}[h!]
  \centering
  \includegraphics[width=0.7\linewidth]{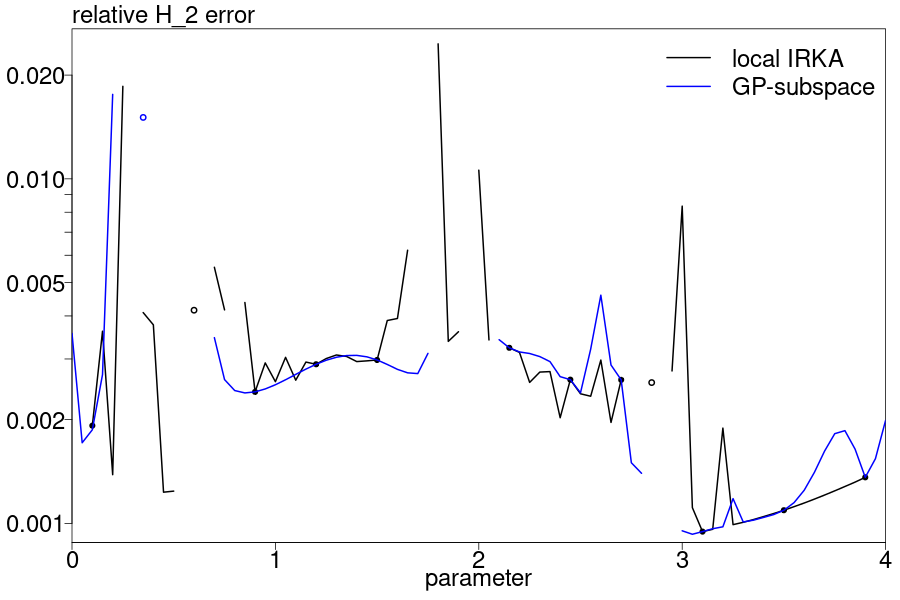}
  \caption{Relative $\mathcal{H}_2$ error for the microthruster. $k=14$.
    Training data are shown as solid points.
    Disconnected test data are shown as hollow points.
  }
  \label{fig:microthruster-k14}
\end{figure}

The microthruster example is just to showcase the accuracy of our proposed method
when combined with a ROM method based on two-sided projection.
The specific combination with IRKA may have several potential issues.
First, IRKA only provides a local optimal ROM,
and there may be an abundance of them depending on the dimensions of the full and the reduced model.
Therefore, different runs of IRKA may give very different pairs of reduced subspaces,
This is reflected in \Cref{fig:microthruster},
as the error curve of local IRKA is occasionally unsmooth.
But for a method that approximates a subspace-valued mapping to work well,
the true mapping needs to be well-defined and smooth in general.
Second, a continuous trajectory of local $\mathcal{H}_2$-optimal ROMs may not be all stable,
which is possible because IRKA may converge to unstable ROMs.
In fact, stability may break multiple times as parameter varies.
Finally, there may not be a continuous trajectory of local $\mathcal{H}_2$-optimal ROMs
across the parameter space, so a good sample of local IRKA subspaces may not exist.
In \Cref{fig:microthruster-k14}, we show some results for $k=14$,
where we use a sample of 10 points for our model.
For the three segments of the parameter space
where the error curve of local IRKA is relatively continuous,
our method is able to maintain the error level,
but overall the error curve is discontinuous and the ROMs can be unstable.
This situation gets worse as $k$ increases in this example.

\section{A limitation of interpolation on tangent space}
\label{sup:discussion}

In general, subspace interpolation is more accurate than the other two interpolation methods.
But when: (1) sample size $l$ is small;
(2) subspace dimension $k$ is large; or (3) parameter dimension $d$ is large,
the accuracy of all these methods can be unsatisfactory.
\cite{Baur2017} Sec. 9.6 also noted that the accuracy of matrix interpolation deteriorates
between sample points when $k$ increases, and gave a tentative explanation.
Here we give an explanation of why interpolation on tangent spaces of a manifold,
which includes subspace and manifold interpolation,
fails in these situations.

When a point $p'$ on a complete Riemannian manifold $\mathcal{M}$ is pulled back
to the tangent space $T_p \mathcal{M}$ of a reference point $p$ via the exponential map,
the preimage $\exp_p^{-1}(p')$ contains an infinite number of tangent vectors.
The Riemannian logarithm $\log_p(p')$ is defined as the smallest tangent vector within this preimage,
which lies in a star-shaped neighborhood of zero called the injectivity domain $\text{ID}(p)$.
When a continuous map $f: \Theta \mapsto \mathcal{M}$ is pulled back to $T_p \mathcal{M}$,
the preimage $(\exp_p^{-1} \circ f)(\Theta)$ may have a connected component in $\text{ID}(p)$,
which can be approximated given enough sample points.
But this component will be increasingly distorted as it approaches the boundary of $\text{ID}(p)$,
called the tangent cut locus $TCL(p)$.
This phenomenon can be observed, for example, in an azimuthal equidistant projection of the Earth.
If the preimage only has connected components that intersects $TCL(p)$ or beyond,
then the map cannot be approximated on $T_p \mathcal{M}$
by continuous maps interpolating points in $\text{ID}(p)$.
As $l$ decreases, $d$ increases, or $k$ increases,
all sample points become further away from each other,
and their Riemannian logarithms move closer to the tangent cut locus for any reference point.
And as $d$ or $k$ increases, the map is more likely to cross the cut locus of any reference point.
Therefore, the map becomes more difficult to approximate on the tangent space in these situations.

\bibliographystyle{mysiamplain-nourl}
\bibliography{gps-main.bib}

\end{document}

%% file: gps-shared.tex
\usepackage{amsmath} %
\usepackage{amssymb} %
\usepackage{amsfonts} %
\usepackage{amsopn} %
\DeclareMathOperator{\diag}{diag}
\usepackage{mathrsfs} %
\usepackage{euscript} %
\usepackage{mathtools}

\newsiamthm{claim}{Claim}
\newsiamremark{remark}{Remark}
\newsiamremark{hypothesis}{Hypothesis}
\crefname{hypothesis}{Hypothesis}{Hypotheses} %
\newsiamremark{assumption}{Assumption}
\crefname{assumption}{Assumption}{Assumptions} %
\usepackage{enumitem}
\setlist[enumerate]{leftmargin=.5in}
\setlist[itemize]{leftmargin=.5in}

\usepackage{algorithm} %
\usepackage{algpseudocode} %
\algnewcommand\algorithmicinput{\textbf{Input:}}
\algnewcommand\algorithmicoutput{\textbf{Output:}}
\algnewcommand\algorithmicnote{\textbf{Note:}}
\algnewcommand\Input{\item[\algorithmicinput]}%
\algnewcommand\Output{\item[\algorithmicoutput]}%
\algnewcommand\Note{\item[\algorithmicnote]}%

\usepackage{threeparttable} %
\usepackage{booktabs} %

\usepackage{graphicx} %
\usepackage{epstopdf}
\ifpdf
  \DeclareGraphicsExtensions{.eps,.pdf,.png,.jpg}
\else
  \DeclareGraphicsExtensions{.eps}
\fi
\usepackage[all]{xy}
\usepackage{tikz-cd}

\usepackage{microtype} %
\usepackage[strict=true]{csquotes} %
\usepackage[strings]{underscore} %

\headers{Gaussian Process Subspace Regression}{R. Zhang, S. Mak, and D. Dunson}

\title{Gaussian Process Subspace Regression for Model Reduction\thanks{Submitted to the editors.
\funding{This work was supported, in part, by the National Science Foundation grant DMS-1638521
and by grant N00014-21-1-2510-01 of the United States Office of Naval Research.}}}

\author{Ruda Zhang\thanks{Department of Mathematics, Duke University, Durham, NC 27710 USA
    (\email{ruda.zhang@duke.edu})}
  \and Simon Mak\thanks{Department of Statistical Science, Duke University, Durham, NC 27710 USA
    (\email{sm769@duke.edu})}
  \and David Dunson\thanks{Department of Mathematics and Department of Statistical Science,
    Duke University, Durham, NC 27710 USA (\email{dunson@duke.edu})}}

%% file: gps-main.bbl
\begin{thebibliography}{10}

\bibitem{Afsari2011}
{\sc B.~Afsari}, {\em Riemannian {$L^p$} center of mass: Existence, uniqueness,
  and convexity}, Proc. Amer. Math. Soc., 139 (2011), pp.~655--673.

\bibitem{Amsallem2010}
{\sc D.~Amsallem}, {\em Interpolation on manifolds of CFD-based fluid and
  finite element-based structural reduced-order models for on-line aeroelastic
  predictions}, phdthesis, Stanford Univ., 2010.

\bibitem{Amsallem2008}
{\sc D.~Amsallem and C.~Farhat}, {\em Interpolation method for adapting
  reduced-order models and application to aeroelasticity}, AIAA J., 46 (2008),
  pp.~1803--1813.

\bibitem{Amsallem2011}
{\sc D.~Amsallem and C.~Farhat}, {\em An online method for interpolating linear
  parametric reduced-order models}, SIAM J. Sci. Comput., 33 (2011),
  pp.~2169--2198.

\bibitem{Antoulas2020}
{\sc A.~C. Antoulas, C.~A. Beattie, and S.~Gugercin}, {\em Interpolatory
  Methods for Model Reduction}, SIAM, 2020.

\bibitem{Baur2017}
{\sc U.~Baur, P.~Benner, B.~Haasdonk, C.~Himpe, I.~Martini, and M.~Ohlberger},
  {\em Chapter 9: Comparison of methods for parametric model order reduction of
  time-dependent problems}, in Model Reduction and Approximation: Theory and
  Algorithms, P.~Benner, M.~Ohlberger, A.~Cohen, and K.~Willcox, eds.,
  Computational Science \& Engineering, SIAM, 2017, pp.~377--407.

\bibitem{Bendokat2020}
{\sc T.~Bendokat, R.~Zimmermann, and P.~A. Absil}, {\em A {Grassmann} manifold
  handbook: Basic geometry and computational aspects}.
\newblock arXiv, 2020.

\bibitem{Benner2015}
{\sc P.~Benner, S.~Gugercin, and K.~Willcox}, {\em A survey of projection-based
  model reduction methods for parametric dynamical systems}, SIAM Rev., 57
  (2015), pp.~483--531.

\bibitem{Chaturantabut2010}
{\sc S.~Chaturantabut and D.~C. Sorensen}, {\em Nonlinear model reduction via
  discrete empirical interpolation}, SIAM J. Sci. Comput., 32 (2010),
  pp.~2737--2764.

\bibitem{chen2020function}
{\sc J.~Chen, S.~Mak, V.~R. Joseph, and C.~Zhang}, {\em Function-on-function
  kriging, with applications to three-dimensional printing of aortic tissues},
  Technometrics,  (2020), pp.~1--12.

\bibitem{Chikuse2003}
{\sc Y.~Chikuse}, {\em Statistics on Special Manifolds}, Springer-Verlag, New
  York, 2003.

\bibitem{Furrer2006}
{\sc R.~Furrer, M.~G. Genton, and D.~Nychka}, {\em Covariance tapering for
  interpolation of large spatial datasets}, J. Comput. Graph. Statist., 15
  (2006), pp.~502--523.

\bibitem{Golub2013}
{\sc G.~H. Golub and C.~F. Van~Loan}, {\em Matrix Computations}, Johns Hopkins
  University Press, Baltimore, Maryland, 2013.

\bibitem{Gramacy2020}
{\sc R.~B. Gramacy}, {\em Surrogates: {Gaussian} Process Modeling, Design, and
  Optimization for the Applied Sciences}, Chapman and Hall/CRC, 2020.

\bibitem{Gramacy2015}
{\sc R.~B. Gramacy and D.~W. Apley}, {\em Local {Gaussian} process
  approximation for large computer experiments}, J. Comput. Graph. Statist., 24
  (2015), pp.~561--578.

\bibitem{Grohs2013}
{\sc P.~Grohs}, {\em Quasi-interpolation in {Riemannian} manifolds}, IMA J.
  Numer. Anal., 33 (2013), pp.~849--874.

\bibitem{Gugercin2008}
{\sc S.~Gugercin, A.~C. Antoulas, and C.~Beattie}, {\em $\mathcal{H}_2$ model
  reduction for large-scale linear dynamical systems}, SIAM J. Matrix Anal.
  Appl., 30 (2008), pp.~609--638.

\bibitem{Guhaniyogi2016}
{\sc R.~Guhaniyogi and D.~B. Dunson}, {\em Compressed {Gaussian} process for
  manifold regression}, J. Mach. Learn. Res., 17 (2016), pp.~1--26.

\bibitem{Halko2011}
{\sc N.~Halko, P.~G. Martinsson, and J.~A. Tropp}, {\em Finding structure with
  randomness: Probabilistic algorithms for constructing approximate matrix
  decompositions}, SIAM Rev., 53 (2011), pp.~217--288.

\bibitem{Higham2005}
{\sc N.~J. Higham}, {\em The scaling and squaring method for the matrix
  exponential revisited}, SIAM J. Matrix Anal. Appl., 26 (2005),
  pp.~1179--1193.

\bibitem{Hokanson2020}
{\sc J.~M. Hokanson and C.~C. Magruder}, {\em $\mathcal{H}_2$-optimal model
  reduction using projected nonlinear least squares}, SIAM J. Sci. Comput., 42
  (2020), pp.~A4017--A4045.

\bibitem{Kaufman2011}
{\sc C.~G. Kaufman, D.~Bingham, S.~Habib, K.~Heitmann, and J.~A. Frieman}, {\em
  Efficient emulators of computer experiments using compactly supported
  correlation functions, with an application to cosmology}, Ann. Appl. Stat., 5
  (2011), pp.~2470--2492.

\bibitem{LinLZ2019}
{\sc L.~Lin, N.~Mu, P.~Cheung, and D.~Dunson}, {\em Extrinsic {Gaussian}
  processes for regression and classification on manifolds}, Bayesian Anal., 14
  (2019), pp.~887--906.

\bibitem{LinLZ2017}
{\sc L.~Lin, B.~S. Thomas, H.~Zhu, and D.~B. Dunson}, {\em Extrinsic local
  regression on manifold-valued data}, J. Am. Stat. Assoc., 112 (2017),
  pp.~1261--1273.

\bibitem{Lumley1967}
{\sc J.~L. Lumley}, {\em The structure of inhomogeneous turbulent flows},
  Atmospheric Turbulence and Radio Wave Propagation,  (1967).

\bibitem{Magnus2019}
{\sc J.~R. Magnus and H.~Neudecker}, {\em Matrix Differential Calculus with
  Applications in Statistics and Econometrics}, Wiley, 2019.

\bibitem{Mak2018}
{\sc S.~Mak, C.-L. Sung, X.~Wang, S.-T. Yeh, Y.-H. Chang, V.~R. Joseph,
  V.~Yang, and C.~F.~J. Wu}, {\em An efficient surrogate model for emulation
  and physics extraction of large eddy simulations}, J. Am. Stat. Assoc., 113
  (2018), pp.~1443--1456.

\bibitem{Mallasto2018}
{\sc A.~Mallasto and A.~Feragen}, {\em Wrapped {Gaussian} process regression on
  {Riemannian} manifolds}, in Proceedings of the IEEE Conference on Computer
  Vision and Pattern Recognition (CVPR), June 2018.

\bibitem{Moore1981}
{\sc B.~C. Moore}, {\em Principal component analysis in linear systems:
  Controllability, observability, and model reduction}, IEEE Trans. Autom.
  Control, 26 (1981), pp.~17--32.

\bibitem{anemometer}
{\sc {MORwiki~Community}}, {\em Anemometer}.
\newblock Model Order Reduction Wiki, 2018.

\bibitem{NiuM2019}
{\sc M.~Niu, P.~Cheung, L.~Lin, Z.~Dai, N.~Lawrence, and D.~Dunson}, {\em
  Intrinsic {Gaussian} processes on complex constrained domains}, J. R. Stat.
  Soc. Ser. B. Stat. Methodol., 81 (2019), pp.~603--627.

\bibitem{thermalmodel}
{\sc {Oberwolfach~Benchmark~Collection}}, {\em Thermal model}.
\newblock Model Order Reduction Wiki, 2018.

\bibitem{Panzer2010}
{\sc H.~Panzer, J.~Mohring, R.~Eid, and B.~Lohmann}, {\em Parametric model
  order reduction by matrix interpolation}, at - Automatisierungstechnik, 58
  (2010), pp.~475--484.

\bibitem{Petersen2019}
{\sc A.~Petersen and H.-G. M{\"{u}}ller}, {\em Fr{\'{e}}chet regression for
  random objects with {Euclidean} predictors}, Ann. Statist., 47 (2019),
  pp.~691--719.

\bibitem{Rasmussen2006}
{\sc C.~E. Rasmussen and C.~K.~I. Williams}, {\em Gaussian Processes for
  Machine Learning}, MIT Press, Cambridge, MA, 2006.

\bibitem{mmess}
{\sc J.~Saak, M.~K{\"{o}}hler, and P.~Benner}, {\em {M-M.E.S.S.} - the matrix
  equation sparse solver library}.
\newblock Zenodo, Apr. 2021.

\bibitem{Sander2016}
{\sc O.~Sander}, {\em Geodesic finite elements of higher order}, IMA J. Numer.
  Anal., 36 (2016), pp.~238--266.

\bibitem{Santner2018}
{\sc T.~J. Santner, B.~J. Williams, and W.~I. Notz}, {\em The Design and
  Analysis of Computer Experiments}, Springer, New York, NY, 2018.

\bibitem{Schmid2008}
{\sc P.~Schmid and J.~Sesterhenn}, {\em Dynamic mode decomposition of numerical
  and experimental data}, in 61st Annual Meeting of the APS Division of Fluid
  Dynamics, vol.~53, 2008.

\bibitem{Son2013}
{\sc N.~T. Son}, {\em A real time procedure for affinely dependent parametric
  model order reduction using interpolation on {Grassmann} manifolds}, Int. J.
  Numer. Methods Eng., 93 (2013), pp.~818--833.

\bibitem{Strang2019}
{\sc G.~Strang}, {\em Linear Algebra and Learning from Data},
  Wellesley-Cambridge Press, 2019.

\bibitem{YangY2016}
{\sc Y.~Yang and D.~B. Dunson}, {\em Bayesian manifold regression}, Ann.
  Statist., 44 (2016), pp.~876--905.

\bibitem{YeK2016}
{\sc K.~Ye and L.-H. Lim}, {\em Schubert varieties and distances between
  subspaces of different dimensions}, SIAM J. Matrix Anal. Appl., 37 (2016),
  pp.~1176--1197.

\bibitem{ZhangRD2020nr}
{\sc R.~Zhang}, {\em Newton retraction as approximate geodesics on
  submanifolds}.
\newblock arXiv, June 2020.

\bibitem{ZhangRD2021nbb}
{\sc R.~Zhang and R.~Ghanem}, {\em Normal-bundle bootstrap}, SIAM J. Math. Data
  Sci., 3 (2021), pp.~573--592.

\bibitem{ZhangRD2020EnvEcon}
{\sc R.~Zhang, P.~Wingo, R.~Duran, K.~Rose, J.~Bauer, and R.~Ghanem}, {\em
  Environmental economics and uncertainty: Review and a machine learning
  outlook}, in Oxford Research Encyclopedia of Environmental Science, Oxford
  University Press, Aug. 2020.

\bibitem{Zimmermann2014}
{\sc R.~Zimmermann}, {\em A locally parametrized reduced-order model for the
  linear frequency domain approach to time-accurate computational fluid
  dynamics}, SIAM J. Sci. Comput., 36 (2014), pp.~B508--B537.

\bibitem{Zimmermann2019}
{\sc R.~Zimmermann}, {\em Manifold interpolation and model reduction}.
\newblock arXiv, 2019.

\bibitem{Zimmermann2020}
{\sc R.~Zimmermann}, {\em Hermite interpolation and data processing errors on
  {Riemannian} matrix manifolds}, SIAM J. Sci. Comput., 42 (2020),
  pp.~A2593--A2619.

\end{thebibliography}
